\documentclass[onefignum,onetabnum]{siamart190516}

\usepackage{amsfonts}
\usepackage{graphicx}
\usepackage{multirow}   
\usepackage{epstopdf}
\usepackage{algorithmic}
\ifpdf
  \DeclareGraphicsExtensions{.eps,.pdf,.png,.jpg}
\else
  \DeclareGraphicsExtensions{.eps}
\fi

\newcommand{\STAB}[1]{\begin{tabular}{@{}c@{}}#1\end{tabular}}  

\newsiamremark{remark}{Remark}
\newsiamremark{hypothesis}{Hypothesis}
\crefname{hypothesis}{Hypothesis}{Hypotheses}
\newsiamthm{claim}{Claim}

\headers{Low-rank Riemannian multigrid line search}{Marco Sutti and Bart Vandereycken}

\title{Riemannian multigrid line search for low-rank problems\thanks{The work of the first author was supported by the SNSF under research project number 163212.}}

\author{Marco Sutti\thanks{Department of Mathematics, University of Geneva, Geneva 1211, Switzerland 
  (\email{Marco.Sutti@unige.ch}, \email{Bart.Vandereycken@unige.ch}).}
  \and Bart Vandereycken\footnotemark[2]}

\usepackage{amsopn}

\ifpdf
\hypersetup{
  pdftitle={Riemannian multigrid line search for low-rank problems},
  pdfauthor={M. Sutti and B. Vandereycken}
}
\fi

\usepackage{mathtools}
\usepackage{enumerate}  
\usepackage{booktabs}

\newcommand{\tr}{^{\mkern-1mu\textsf {T}}}    
\newcommand{\kt}{\ast\tr}  
\newcommand{\R}{\mathbb{R}}
\newcommand{\dx}{\,\mathrm{d}x}
\newcommand{\dy}{\,\mathrm{d}y}

\newcommand{\wex}{w_{\mathrm{ex}}}
\newcommand{\epsmach}{\varepsilon_{\mathrm{mach}}}

\DeclareMathOperator{\grad}{grad}
\DeclareMathOperator{\D}{D}   
\renewcommand*{\P}{\operatorname{P}}   

\newcommand{\cM}{{\mathcal M}}
\newcommand{\cMk}{\cM_{k}} 
\newcommand{\cMkh}{\cM_{h}^{k}} 
\newcommand{\cMkH}{\cM_{H}^{k}} 
\newcommand{\cMh}{\cM_{h}} 
\newcommand{\cMH}{\cM_{H}} 
\newcommand{\cI}{{\mathcal I}}
\newcommand{\cF}{{\mathcal F}}
\newcommand{\F}{\mathrm{F}}

\DeclareMathOperator{\rank}{rank}
\DeclareMathOperator{\diag}{diag}
\DeclareMathOperator{\blkdiag}{blkdiag}
\DeclareMathOperator*{\trace}{\mathrm{tr}}

\begin{document}

\maketitle

\begin{abstract} 
Large-scale optimization problems arising from the discretization of problems involving PDEs sometimes admit solutions that can be well approximated by low-rank matrices. In this paper, we will exploit this low-rank approximation property by solving the optimization problem directly over the set of low-rank matrices. In particular, we introduce a new multilevel algorithm, where the optimization variable is constrained to the Riemannian manifold of fixed-rank matrices. In contrast to most other multilevel algorithms where the rank is chosen adaptively on each level in order to control the perturbation due to the low-rank truncation, we can keep the ranks (and thus the computational complexity) fixed throughout the iterations. Furthermore, classical implementations of line searches based on Wolfe conditions enable computing a solution where the numerical accuracy is limited to about the square root of the machine epsilon. Here, we propose an extension to Riemannian manifolds of the line search of Hager and Zhang, which uses approximate Wolfe conditions that enable computing a solution on the order of the machine epsilon. Numerical experiments demonstrate the computational efficiency of the proposed framework.
\end{abstract} 

\begin{keywords}
low-rank matrices, optimization on manifolds,  multilevel optimization, Riemannian manifolds, retraction-based optimization, line search, roundoff error
\end{keywords}

\begin{AMS}
  65F10, 65N22, 65F50, 65K10
 \end{AMS}

\section{Introduction}

The topic of this paper is the efficient solution of certain large-scale variational problems arising from the discretization of elliptic PDEs. In particular, we combine Riemannian optimization on the manifold of fixed-rank matrices with ideas from nonlinear multigrid and multilevel optimization. The low-rank manifold will allow us to approximate the solution with significantly less degrees of freedom. In addition, the idea of recursive coarse-grid corrections from multigrid will lead to almost mesh-independent convergence of our algorithm similar to classical multigrid algorithms.

Approximating very large matrices by low rank is a popular technique to speed up numerical calculations. In the context of high-dimensional problems, this is done in so-called low-rank matrix and tensor methods, where tensors are the higher order analog of two-dimensional matrices~\cite{Hackbusch:2012}. One of the early examples are low-rank solvers for the Lyapunov equation, $AX+XA\tr = C$, and other matrix equations; see~\cite{Simoncini:2016} for a recent overview. In order to obtain a low-rank approximation of the unknown solution $X$, an iterative method has to be used that directly constructs the low-rank approximation. Of particular importance for this paper are methods that accomplish this via Riemannian optimization~\cite{AMS:2008}: the minimization problem (obtained after a possible reformulation of the original problem) is restricted to the manifold of fixed-rank matrices, thereby guaranteeing a low-rank representation of critical points. Examples of such methods are~\cite{Mishra:2011b,Vandereycken:2013,Steinlechner:2016aa} for matrix and tensor completion, \cite{Shalit:2012} for metric learning, \cite{VandereyckenV_2010,Mishra_V_2014,Kressner:2016} for matrix and tensor equations, and \cite{Rakhuba:2018ab,RAKHUBA2019718} for eigenvalue problems. In the context of discretized PDEs these optimization problems are very ill-conditioned, making simple first-order methods like gradient descent unmanageably slow. In \cite{VandereyckenV_2010,Kressner:2016,Rakhuba:2018ab}, for example, the gradient is therefore preconditioned with the inverse of the local Hessian. Solving these Hessian equations is done by a preconditioned iterative scheme, thereby mimicking the class of quasi or truncated Newton methods. We also refer to~\cite{UschmajewV:2019} for a recent overview of geometric methods for obtaining low-rank approximations.

Multilevel optimization is the extension of multigrid, and in particular, the full approximation scheme (FAS) to unconstrained optimization. The MG/Opt method from~\protect{\cite{Nash:2000,Lewis:2005} introduced the idea how to modify} the objective functions on each scale so that they correspond to FAS coarse-grid corrections. Several extensions and theoretical convergence proofs were proposed, including optimization with trust-regions~\cite{Toint:2009} and line searches~\cite{Wen:2009}. 
Related to this paper is the low-rank multigrid method from~\cite{Grasedyck:2007} for matrix equations arising from the discretization of elliptic PDEs. It applies a low-rank approximation after every step of the classical multigrid algorithm from~\cite{Penzl:1997} for the linear Sylvester matrix equation. 
A similar multigrid approach with truncation of the matrix iterates to low rank is used by~\cite{Elman:2018} for the solution of large linear systems of equations arising from the finite element discretization of stochastic PDEs.

Our proposed method is different in the sense that it is closer to MG/Opt and other multilevel optimization algorithms in that it works directly with the manifold of fixed-rank matrices. All the classical components of multigrid are present, plus the additional components from Riemannian optimization that allow us to cope with the curvature of the manifold, and thus to generalize the existing Euclidean algorithm to manifolds. Moreover, our implementation permits us to keep the ranks (and thus the computational complexity) fixed throughout the iterations.

This paper is structured as follows.  We first recall important ideas from multilevel optimization and the geometry of fixed-rank matrices that will be needed later on. The main contribution is in \cref{sec:multilevel_opt_Riem} where we present our new algorithm entitled Riemannian multigrid line search (RMGLS). The presentation will be sufficiently general to be applicable to any multilevel hierarchy of manifolds but the implementation will be explained only for low-rank matrices. In \cref{sec:linesearch}, we discuss the numerical difficulty and our solution to obtain critical points with high accuracy using only first-order information in standard line-search methods. Numerical experiments for both a linear and a nonlinear variational problem are presented in \cref{sec:variational_problems}. Finally, in \cref{sec:comparison}, we compare our method to other low-rank and multilevel methods.

\section{Preliminaries on multilevel optimization and geometry of fixed-rank matrices}

As mentioned above, our algorithm is a generalization of known (Euclidean) multilevel algorithms to Riemannian manifolds. It will then be able to calculate low-rank approximations for the variational problems discussed in \cref{sec:variational_problems} by minimizing a cost function over the manifold of fixed-rank matrices. Before we present this algorithm in~\cref{sec:multilevel_opt_Riem}, we briefly recall two important concepts for its derivation: MG/Opt~\cite{Nash:2000}, a variant of multigrid for optimization problems, and retraction-based Riemannian optimization~\cite{AMS:2008}, a local optimization method well suited to minimize over the set of fixed-rank matrices.

\subsection{Multilevel optimization in Euclidean space}\label{ref:sec_Euclidean_multilevel} 

In this paper, we assume some basic knowledge of (linear) multigrid and the full approximation scheme (FAS). For an introduction to these methods, we refer the reader to the books \protect{\cite{Hackbusch:1985,Trottenberg:2000}}.
While multigrid methods were originally devised to solve large-scale linear systems arising from the discretization of PDEs, they can also be used to solve certain nonlinear problems; see \protect{\cite[Chap.~9]{Hackbusch:1985}} and \protect{\cite[Chap.~5.3]{Trottenberg:2000}}. Indeed, the two fundamental multigrid principles, error smoothing and coarse-grid correction, are also present in FAS when solving nonlinear problems. Error smoothing properties are typically encountered in elliptic PDEs and other problems that show some degree of ellipticity. The main difference with respect to linear multigrid is that FAS does not use an error equation on the coarse grid but aims to find a correction based on the full unknown.
The residual equation of FAS on a fine grid is given by~\protect{\cite[eq.~(5.3.11)]{Trottenberg:2000}}
\begin{equation}\label{eq:FAS_defect_fine_grid}
   N_{h} ( w_{h} ) = r_{h} + N_{h} ( \bar{u}_{h} ),
\end{equation}
where $ N_{h} $ is a discrete nonlinear operator, $ w_{h}= \bar{u}_{h} + e_{h} $ is the full approximation, $ \bar{u}_{h} $ is the smoothed approximation, and $ e_{h} $ is the correction to be computed.
On the coarse grid, equation~\eqref{eq:FAS_defect_fine_grid} is approximated by~\protect{\cite[eq.~(5.3.12)]{Trottenberg:2000}}
\[
   N_{H} ( w_{H} ) = r_{H} + N_{H} (\bar{u}_{H}),
\]
with $ w_{H} = \bar{u}_{H} + e_{H} $ the full approximation on the coarse grid, and $ \bar{u}_{H} $ the restriction of $ \bar{u}_{h} $ to the coarse grid. This equation has to be solved for the coarse-grid correction $ e_{H} $, which is computed as the difference of $ \bar{u}_{H} $ and $ w_{h} $, and is then transferred to the fine grid as a correction $ e_{h} $.
In FAS, smoothing is achieved via a nonlinear relaxation procedure having appropriate error smoothing properties, such as nonlinear Gauss--Seidel or weighted Jacobi.

FAS can be generalized to a multilevel algorithm for minimizing a differentiable objective function $ f $. The original idea goes back to the MG/Opt algorithm \protect{\cite{Nash:2000,Lewis:2005}}. Here, we briefly explain the main idea for two grids since the algorithm on more grids is recursively defined from it, and we will explain the algorithm for Riemannian manifolds in more detail in~\cref{sec:multilevel_opt_Riem}. Let the subscripts $ \cdot_{h} $, $ \cdot_{H} $ denote quantities on the fine $\Omega_h \simeq \R^n$ and the coarse grid $\Omega_H\simeq \R^N$ , respectively. Let $ f_{h} \colon \Omega_h \to \R$ be our original  objective function $f$ that we optimize with an initial guess $ \bar{x}_{h} \in \Omega_h$ that is sufficiently smoothed. As in FAS, we introduce a modification to the coarse-grid objective function $ f_{H} $.
Let $ g^{\mathrm{E}}(z_{1},z_{2}) \coloneqq z_{1}\tr z_{2} $ denote the Euclidean inner product and $ I^{H}_{h} \colon \Omega_h \to \Omega_H$ the restriction operator. At iteration $ i $ of MG/Opt, let $ {x}_{H}^{(i)} = I_{h}^{H} \bar{x}_h \in \Omega_H$ be the iterate on the coarse grid.
Then by minimizing the model
\begin{equation}\label{eq:psi_H}
   \psi_{H}  \colon \Omega_H \to \R, \qquad x_H \mapsto f_{H}(x_{H}) - g^{\mathrm{E}}(x_{H}, \kappa_{H}),
\end{equation}
with
\begin{equation}\label{eq:kappa_H_cgc}
   \kappa_{H} \coloneqq \nabla f_{H}({x}_{H}^{(i)}) - I_{h}^{H} \nabla f_{h}(\bar{x}_{h}), 
\end{equation}
one obtains a two-grid cycle for optimizing $ f_{h} $.
On the coarser level, the minimization of \cref{eq:psi_H} starts at the smoothed approximation $ {x}_{H}^{(i)} $, hence we can rewrite \cref{eq:psi_H} in the following way: find an update $ e_{H} $ such that
\begin{equation}\label{eq:psi_H_CGC}
   \psi_{H}({x}_{H}^{(i)}+e_{H}) \coloneqq f_{H}({x}_{H}^{(i)}+e_{H}) - g^{\mathrm{E}}({x}_{H}^{(i)}+e_{H}, \kappa_{H} )
\end{equation}
is sufficiently minimized at ${x}_{H}^{(i+1)} = {x}_{H}^{(i)}+e_{H}$. This coarse-grid update $e_H$ is then transported back to the fine grid using the interpolation operator $ I^{h}_{H} \colon \Omega_H \to \Omega_h$, and used to correct $\bar{x}_h$.

The linear modification~\eqref{eq:psi_H} to $ f_{H} $ is one of the central tenets of multilevel optimization, as proposed in the MG/Opt method of \protect{\cite{Nash:2000,Lewis:2005}} and similar multilevel algorithms in~\cite{Gratton:2008,Wen:2009}. The model $ \psi_{H} $ is actually a generalization of the coarse-grid correction equation of the FAS scheme in the context of optimization. Indeed, applying FAS for solving the nonlinear critical point equation $\nabla f_h(x) = 0$ at the approximation ${x}_{h}^{(i)}$ gives the coarse-grid correction \protect{\cite[Chap.~(5.3.4)]{Trottenberg:2000}}
\[
   \nabla f_H ( {x}_{H}^{(i)} + e_{H} )  -  \nabla f_H ({x}_{H}^{(i)}) = -I_h^H \nabla f_h (\bar{x}_{h})
\]
that has to be solved for $e_H$. A solution of this equation can be trivially written as
\[
   \nabla f_H ( {x}_{H}^{(i)} + e_{H} ) -  (\ \nabla f_H ({x}_{H}^{(i)}) -I_h^H \nabla f_h (\bar{x}_{h}) \ )= 0,
\]
which is exactly a critical point of~\cref{eq:psi_H_CGC} with definition~\cref{eq:kappa_H_cgc} for $\kappa_H$.

As in classical multigrid methods, the error has to be smooth in order to be representable on the coarse grid. For classical multigrid or FAS, iterative methods such as weighted Jacobi and Gauss--Seidel, and their nonlinear versions, can be used to smooth the error. Analogously, in the optimization framework, one can use cheap first-order optimization methods. Practice has shown that weighted versions of steepest descent, coordinate search and limited memory BFGS are effective smoothers for a wide range of large-scale multilevel optimization problems; see, e.g., \protect{\cite{Gratton:2010}}.

\begin{figure}[htbp]
   \centering
   \includegraphics[width=0.75\columnwidth]{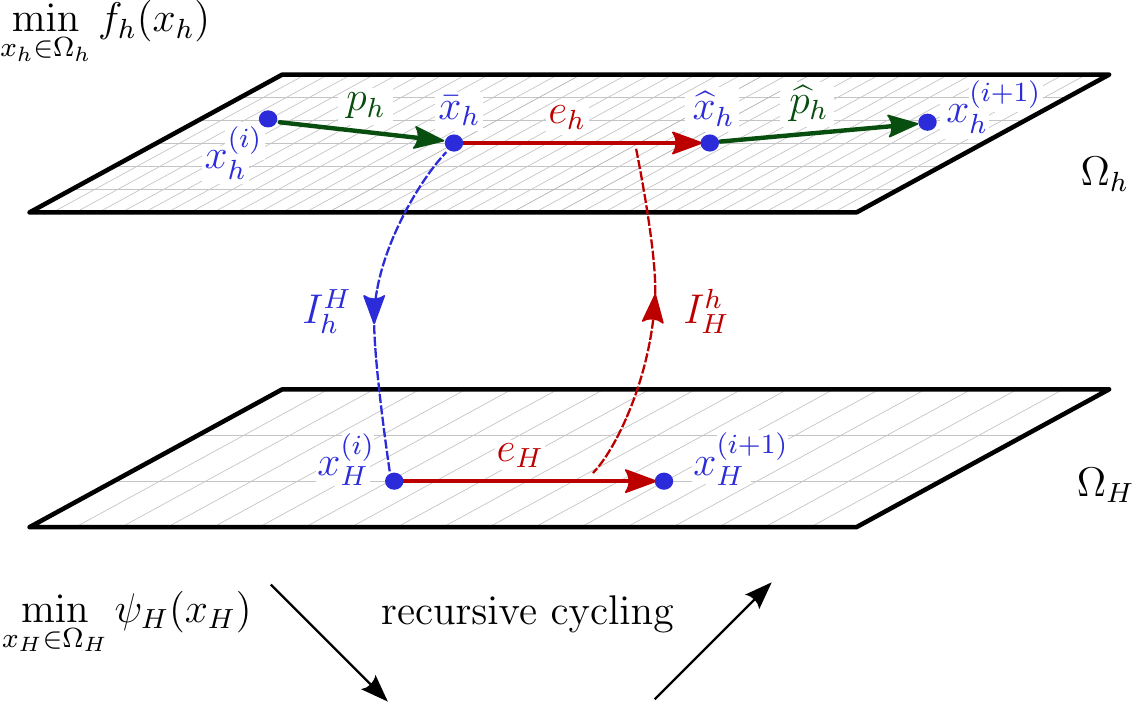}
   \caption{A two-grid cycle for minimizing an objective function.}\label{fig:scheme_of_MGOPT}
\end{figure}

Except for the introduction of the model \cref{eq:psi_H}, the principle behind the multigrid two-grid cycle remains the same in the optimization context. \Cref{fig:scheme_of_MGOPT} illustrates the two-grid cycle of a multilevel optimization scheme. The initial guess at iterate $ i $ is denoted by $ x_{h}^{(i)} $ and the pre-smoothing update by $ p_{h} $, likewise $ \widehat{p}_{h} $ is the post-smoothing update, resulting in the next iterate $ x_{h}^{(i+1)} $.
In the next section, we will generalize this two-grid optimization cycle (and figure) to Riemannian manifolds.

\subsection{The manifold of fixed-rank matrices}

Computing a rank-$k$ approximation to a matrix $ X \in \R^{m \times n} $ can be seen as an optimization problem on the manifold of fixed-rank matrices
\[
    \cMk = \lbrace X \in \R^{m \times n} \colon \rank(X) = k \rbrace.
\]
Using the SVD, one has the equivalent characterization 
\begin{equation*}
    \begin{split}
       \cMk = \lbrace U\Sigma V\tr \colon & U \in \mathrm{St}^{m}_{k}, \ V \in \mathrm{St}^{n}_{k}, \\
              & \Sigma = \diag(\sigma_{1}, \sigma_{2}, \ldots, \sigma_{k} ) \in \R^{k \times k}, \ \sigma_{1} \geq \cdots \geq \sigma_{k} > 0 \rbrace,
    \end{split}
\end{equation*}
where $ \mathrm{St}^{m}_{k} $ is the Stiefel manifold of $ m \times k $ real matrices with orthonormal columns, and  $ \diag(\sigma_{1}, \sigma_{2}, \ldots, \sigma_{k} ) $ is a square matrix with $ \sigma_{1}, \sigma_{2}, \ldots, \sigma_{k} $ on its diagonal.
The following proposition shows that $ \cMk $ is indeed a smooth manifold and has a compact representation for its tangent space.

\begin{proposition}[\protect{\cite[Prop.~2.1]{Vandereycken:2013}}]
   The set $ \cMk $ is a smooth submanifold of dimension $ (m+n-k)k $ embedded in $ \R^{m \times n} $. Its tangent space $ T_{X}\cMk $ at $ X = U\Sigma V\tr \in \cMk $ is given by
   \begin{equation}\label{eq:tan_vec_format}
      T_{X}\cMk =
      \begin{bmatrix}
         U  &  U_{\perp}
      \end{bmatrix}
      \begin{bmatrix}
         \R^{k \times k}      &  \R^{k \times (n-k)}  \\
         \R^{(m-k) \times k}  &  0_{(m-k)\times (n-k)}
      \end{bmatrix}
      \begin{bmatrix}
         V  &  V_{\perp}
      \end{bmatrix}\tr.  
   \end{equation}
   In addition, every tangent vector $\xi \in T_{X}\cMk$ can be written as
\begin{equation}\label{eq:tan_vec_format_small_param}
\xi = UMV\tr + U_{\mathrm{p}}V\tr + UV_{\mathrm{p}}\tr,
\end{equation}
with $M \in \R^{k\times k}$, $U_{\mathrm{p}} \in \R^{m\times k}$, $V_{\mathrm{p}} \in \R^{n \times k}$ such that $U_{\mathrm{p}}\tr U = V_{\mathrm{p}}\tr V = 0$.
\end{proposition}
Observe that since $ \cMk \subset \R^{m \times n} $, we represent tangent vectors in~\eqref{eq:tan_vec_format} and~\eqref{eq:tan_vec_format_small_param} as matrices of the same dimensions.
The Riemannian metric is the restriction of the Euclidean metric on $ \R^{m \times n} $ to the submanifold $ \cMk $,
\[
    g_{X}(\xi,\eta) = \langle \xi, \eta \rangle = \trace(\xi\tr \eta), \quad \text{with} \ X \in \cMk \ \text{and} \ \xi, \eta \in T_{X}\cMk.
\]
The Riemannian gradient of a smooth function $ f \colon \cMk \to \R $ at $ X \in \cMk $ is defined as the unique tangent vector $ \grad f(X) $ in $ T_{X}\cMk $ such that
\[
   \langle \, \grad f(X), \xi \, \rangle = \D f(X) [\xi] \quad \text{for all} \ \xi \in T_{X}\cMk,
\]
where $ \D f $ denotes the directional derivatives of $ f $.
More concretely, for embedded submanifolds, the Riemannian gradient is given by the orthogonal projection onto the tangent space of the Euclidean gradient of $ f $ seen as a function on the embedding space $ \R^{m \times n} $; see, e.g., \protect{\cite[eq.~(3.37)]{AMS:2008}}.
Defining $ \P_{U} = UU\tr $ and $ \P_{U}^{\perp} = I - \P_{U} $ for any $ U \in \mathrm{St}^{m}_{k} $, the orthogonal projection onto the tangent space at $ X $ is \protect{\cite[eq.~(2.5)]{Vandereycken:2013}}
\[
    \P_{T_{X}\cMk} \colon \R^{m \times n} \to T_{X}\cMk, \quad Z \mapsto \P_{U} Z \P_{V} + \P_{U}^{\perp} Z \P_{V} + \P_{U} Z \P_{V}^{\perp}.
\]
Then, denoting $ \nabla f(X) $ the Euclidean gradient of $ f $ at $ X $, the Riemannian gradient is given by
\begin{equation}\label{eq:Riemannian_grad_as_projection}
    \grad f(X) = \P_{T_{X}\cMk} \! \big( \nabla f(X) \big).
\end{equation}

\subsection{The orthographic retraction} \label{sec:ortho_retr}
A retraction is a smooth map from the tangent space to the manifold, $ R_{X}\colon T_{X}\cMk \to \cMk $, used to map tangent vectors to points on the manifold. It is, essentially, any smooth first-order approximation of the exponential map of the manifold; see, e.g., \protect{\cite[Definition~1]{Absil:2012}}. In order to establish convergence of the Riemannian algorithms, it is sufficient for the retraction to be defined only locally.

\begin{figure}[htbp]
   \centering
   \includegraphics[width=0.5\columnwidth]{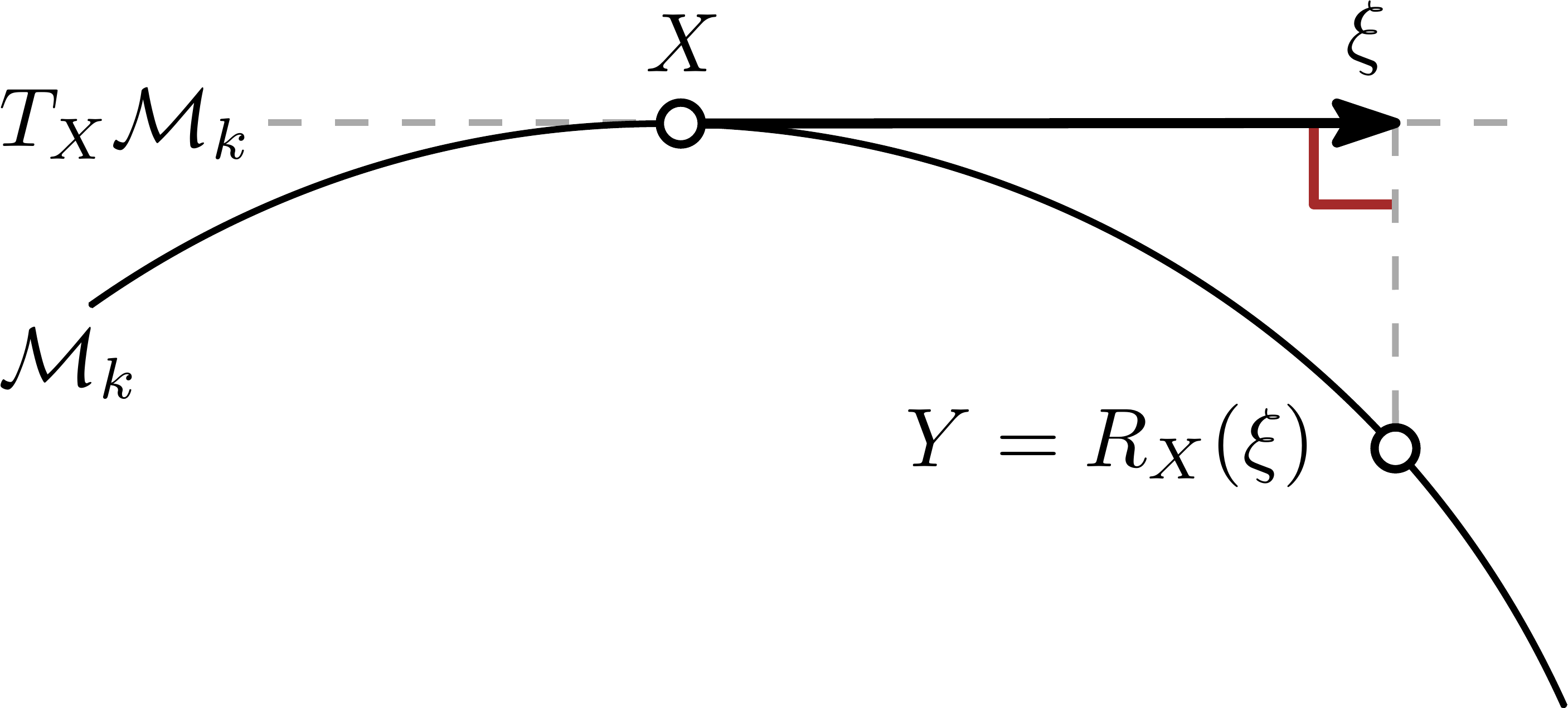}
   \caption{The orthographic retraction.}\label{fig:orthographic_retraction}
\end{figure}

In our setting, we have chosen the orthographic\footnote{The name is due to the fact that on a sphere it relates to the orthographic projection known in cartography \protect{\cite[\S4.1]{Absil:2012}}. As a retraction, it goes back as early as~\protect{\cite{Rosen:1961,Luenberger:1972}}.} retraction on $ \cMk $; see~\protect{\cite[\S4.1]{Absil:2012}} and \protect{\cite[\S3.2]{Absil:2015}}. It is defined by setting $ R_{X}(\xi) $ as the point
nearest to $ X + \xi $ in
\[
   X + \xi + T_{X}^{\perp}\cMk \cap \cMk.
\]
If $ \xi $ is sufficiently small, then $ R_{X}(\xi) $ is unique.
Given a point $ X = U \Sigma V\tr \in \cMk $ and a tangent vector $ \xi $ in the format \cref{eq:tan_vec_format_small_param}, the retraction of $ \xi $ at $ X $ is given by \protect{\cite[§3.2]{Absil:2015}}
\begin{equation}\label{eq:ortho_graphic_R}
    R_{X}(\xi) = [U(\Sigma+M)+U_{\mathrm{p}}]\ (\Sigma+M)^{-1} \ [(\Sigma+M)V\tr + V_{\mathrm{p}}\tr].
\end{equation}
\Cref{fig:orthographic_retraction} illustrates the orthographic retraction. 
As a special case, observe that if $X$ is full rank, then $ UU\tr = U\tr U = I $ and $ VV\tr = V\tr V = I $, therefore $ U_{\mathrm{p}} = 0 $ and $ V_{\mathrm{p}} = 0 $, so $R_{X}(\xi) = U (\Sigma+M) V\tr = X + \xi$.

The reason for this choice is that the orthographic retraction has an explicit expression for its inverse. In particular, it satisfies
\begin{equation}\label{eq:inv_ortho_retr}
   \xi \coloneqq R^{-1}_{X}(Y) = \P_{T_{X}\cMk}(Y-X) = \P_{T_{X}\cMk}(Y) - X.
\end{equation}
Equivalently, this can be written in tangent vector format \cref{eq:tan_vec_format_small_param} with the factors
\[
   M_{\xi} = U\tr Y V - \Sigma, \qquad U_{\mathrm{p},\xi} = (I - U U\tr) YV, \qquad V_{\mathrm{p},\xi} = (I - VV\tr) Y\tr U.
\]

When implementing $ R_{X} $ and $ R_{X}^{-1} $, it is important to exploit the factored forms of the rank-$k$ matrices $X$ and $Y$, and the parametrization~\cref{eq:tan_vec_format_small_param} of the tangent vector $\xi$. In that case, the flop counts of $ R_{X} $ and $ R^{-1}_{X} $ are both $ O(nk^{2} + k^{3}) $. See also \protect{\cite[§3.2]{Absil:2015}}.

\section{Riemannian multigrid  line search for low-rank matrices}\label{sec:multilevel_opt_Riem} In this section, we describe the central contribution of our paper: a Riemannian multilevel line-search algorithm, called RMGLS, for the approximate low-rank solution of optimization problems. We detail how the two-grid optimization cycle of MG/Opt can be generalized to the retraction-based framework for the geometry of fixed-rank matrices, both of which were described in the previous section. 

Our algorithm involves the classical components of multigrid (smoothers, prolongation and restriction operators, and a coarse-grid correction) and Riemannian optimization (line search, retractions, gradients). Since this generalization is possible for other types of manifolds, we have presented it with general manifolds in mind. However, remarks on the implementation apply only to the manifold of fixed-rank matrices.

\subsection{Description of the scheme}

\begin{figure}[htbp]
\centering
\includegraphics[width=0.80\columnwidth]{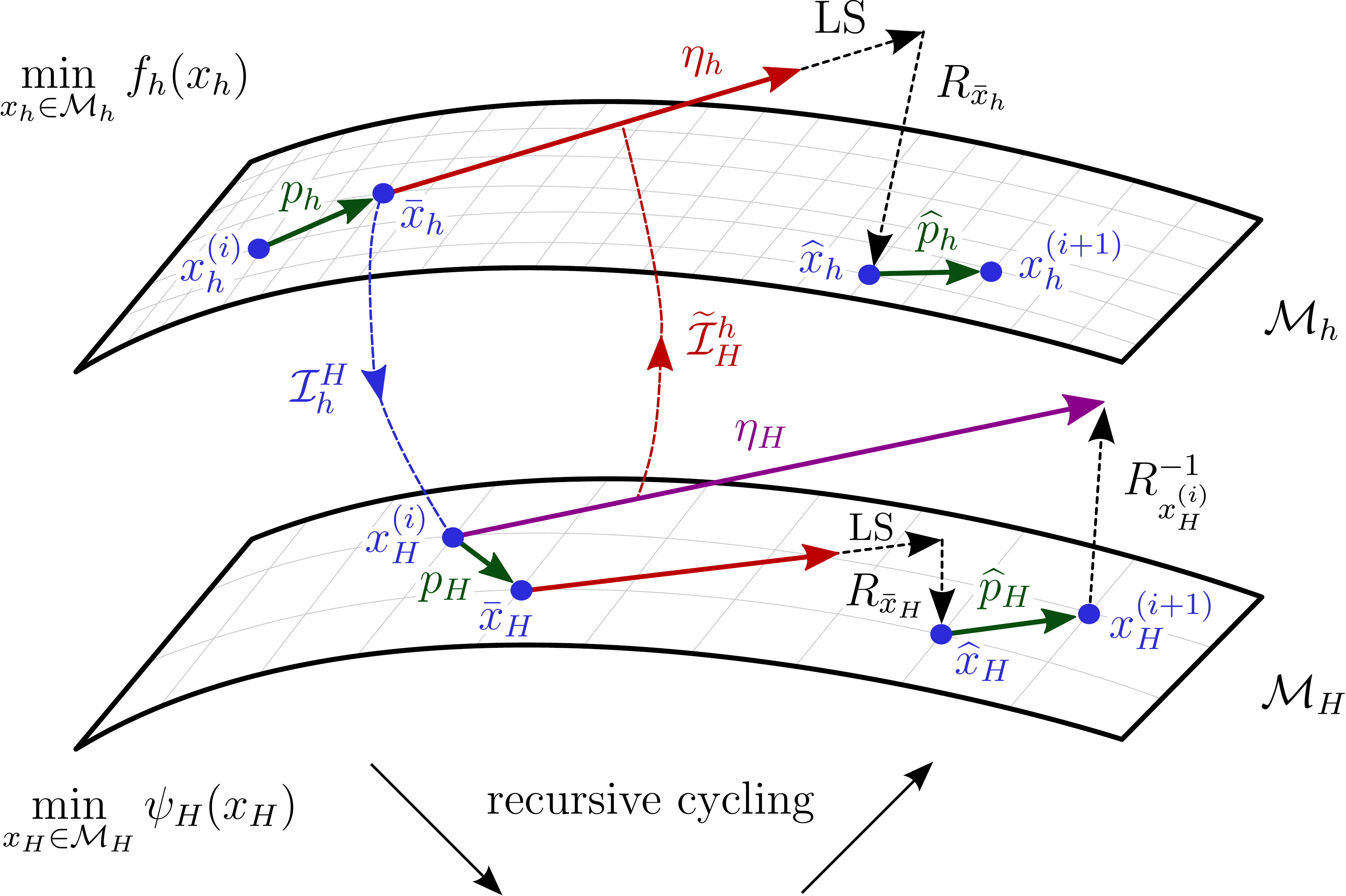}
\caption{The Riemannian multigrid line-search (RMGLS) scheme. The coarse-grid correction is computed either directly or by a recursive application of RMGLS. It is instructive to compare this figure to the Euclidean version in~\cref{fig:scheme_of_MGOPT}.}\label{fig:scheme_of_RMGLS}
\end{figure}

We first describe the algorithm for a two-grid cycle, making reference to \cref{fig:scheme_of_RMGLS}. Recall that quantities related to the fine grid and to the coarse grid are denoted by the subscripts $ \cdot_{h} $ and $ \cdot_{H} $, respectively. For example, $ \cMh $ and $ \cMH $ are the fine and coarse-scale manifolds, respectively. In the experiments presented in section~\ref{sec:variational_problems}, 
we consider a two-dimensional square grid $ \Omega_{h} $ at the finest level $ \ell_{\mathrm{f}} $ with uniform mesh width $ h = 2^{-\ell_{\mathrm{f}}}$, and $ n = 2^{\ell_{\mathrm{f}}} + 1 $ grid points in each direction. A solution on this grid can then be represented as an $ n $-by-$ n $ matrix and its rank-$k$ approximation represents an element of the fine-scale manifold $ \cMh $. Similarly for the coarse-scale manifold $ \cMH $, with $ H = 2h $.

Starting from an approximation $ x_{h}^{(i)} $ on $ \cMh $, we first perform some pre-smoothing steps, then the smoothed approximation $ \bar{x}_{h} $ is restricted to $ \cMH $. This gives us $ x_{H}^{(i)} $, for which we compute a correction $\eta_H$ on $ \cMH $. If the dimension of $ \cMH $ is sufficiently small, $\eta_H$ is computed directly with a trust-region method to minimize $\psi_H$. Otherwise, it is the inverse retraction of the result $x^{(i+1)}_H$ obtained from the recursive application of the two-grid scheme with $ \cMH $ as fine-scale manifold. In the figure, the latter option is depicted for illustration, including the steps performed on $ \cMH $. In both cases, the interpolation $\eta_h$ of the coarse-scale correction $\eta_H$ to the fine scale is applied to $\bar{x}_h$ via line search. The updated approximation $\widehat x_h$ is then post-smoothed and we finally obtain $ x_{h}^{(i+1)} $ as result of one iteration of RMGLS.

An important difference compared to multilevel optimization on Euclidean space is the explicit difference between the approximations $ x_h^{(i)} $, $ \bar{x}_h $, $ \widehat{x}_h $, $ x_h^{(i+1)} $, $ x^{(i)}_H $, $ x^{(i+1)}_H $ that are points on the manifolds $\cMh$ and $\cMH$, and the updates and corrections $ p_h $, $ \widehat{p}_h $, $ \eta_h $, $ \eta_H $ that are tangent vectors on the tangent spaces of $\cMh$ and $\cMH$. This is also clearly visible in \cref{fig:scheme_of_RMGLS} where the approximations are depicted as full circles and tangent vectors as arrows.

In the next subsections, we will explain every component of the algorithm, except for the line search, which will be explained in more detail in \cref{sec:linesearch}. The final algorithm in pseudo-code is listed in~\cref{sec:RMGLS_code}.

\subsection{Tensor-product multigrid} \label{sec:tensor_product_mg}

Observe that a matrix in $ \R^{n \times n} $ can be regarded as an element of the tensor-product space $ \R^{n} \otimes \R^{n} \simeq \R^{n \times n} $. Starting from this observation, it is possible to construct a multigrid algorithm by taking tensor products of standard multigrid components. This approach is known as \emph{tensor-product multigrid} \protect{\cite{Rosen:1995,Penzl:1997}}.

For example, let $ I_{h}^{H} \colon \R^{n} \to \R^{N} $ and $ I_{H}^{h} \colon \R^{N} \to \R^{n} $ denote the standard restriction and prolongation operators for a linear multigrid algorithm with $\R^n$ the fine and $\R^N$ the coarse grid. Let $ \ell_{\mathrm{f}} $ denote the fine level, $ h = 2^{-\ell_{\mathrm{f}}} $, $ H = 2h $, $ n = 2^{\ell_{\mathrm{f}}} - 1 $, and $ N = 2^{\ell_{\mathrm{f}}-1} - 1 $. Then in 1D the restriction $ I_{h}^{H} $ is the $ N \times n $ injection matrix
\begin{equation*}
   (I_{h}^{H})_{ij} =
   \begin{cases}
      1,            & \text{if } j = 2i; \\
      0,            & \text{otherwise}.
   \end{cases}
\end{equation*}
Higher-order extensions for $I_{h}^{H}$ and $I_{H}^{h}$, like full weighting and linear interpolation, are defined analogously; see~\cite{Trottenberg:2000}. Following the tensor-product idea, we can then easily construct a \emph{restriction operator} on the space of matrices by applying $I_{h}^{H}$ to the rows and columns of $X$, 
\begin{equation}\label{eq:restrict_tensor}   \cI_{h}^{H} \colon \R^{n\times n} \to \R^{N \times N}, \quad X \mapsto I_{h}^{H} X (I_{h}^{H})\tr.
\end{equation}
Likewise, an \emph{interpolation operator} for matrices is constructed as
\begin{equation*}
   \cI_{H}^{h} \colon \R^{N \times N} \to \R^{n \times n}, \quad X \mapsto I_{H}^{h} X (I_{H}^{h})\tr.
\end{equation*}
Hence, we have obtained transfer operators between the fine and coarse grids $\R^{n \times n}$ and $\R^{N \times N}$, respectively.

\subsection{Riemannian transfer operators} \label{sec:transfer_operators} 
In our setting, the transfer operators from above are to be applied to rank-$k$ matrices. Let us denote these manifolds by $ \cMkh \subset \R^{n \times n}$ and $ \cMkH \subset \R^{N \times N}$. 

First, we can directly compute the restriction from $ \cMkh $ to $ \cMkH $ by~\eqref{eq:restrict_tensor} since both manifolds are embedded in matrix space. It is clear from~\eqref{eq:restrict_tensor} that $\rank (\cI_{h}^{H}(X_h)) \leq k$ if $X_h$ is a rank-$k$ matrix. In numerical calculations, the rank of $\cI_{h}^{H}(X_h)$ is always equal to $k$, but if it were strictly less we could simply reduce the defining rank of the coarse manifold.\footnote{In the next step of our algorithm RMGLS, the rank of the coarse iterate will typically grow after smoothing and we can then again continue with $ \cMkH $ as our coarse manifold.} The computation of  $ \cI_{h}^{H}(X_{h}) $ is carried out directly on its factorized SVD form, and followed by a reorthogonalization to preserve the SVD format of the result. The entire procedure is summarized in the following box.

\medskip
\fbox{\begin{minipage}{0.90\columnwidth}
\textbf{Restriction} of $X_h = U_h \Sigma_h V_h\tr \in \cMkh$:
\begin{enumerate}[\qquad(1)]
\item Compute \textbf{compact QRs}: \quad $ Q_{U}R_{U} = I_{h}^{H} U_{h} $ and $ Q_{V}R_{V} = I_{h}^{H} V_{h} $

\item Compute \textbf{compact SVD}: \quad $ \widehat{U}\widehat{\Sigma}\widehat{V}\tr = R_{U} \Sigma_{h} R_{V}\tr $

\item Compute \textbf{factors}: \quad $ U_{H} = Q_{U} \widehat{U}, \ \Sigma_{H} = \widehat{\Sigma}, \ V_{H} = Q_{V} \widehat{V} $
\end{enumerate}
Result is $X_{H} = U_H \Sigma_H V_H\tr \in \cMH^{\bar k}$ in SVD form, with $\bar k = \rank ( X_{H} ) $.
\end{minipage}}
\medskip

Next, when transferring tangent vectors between manifolds of different scales, the result of the transfer operators is not necessarily in the tangent space at the transferred points. We therefore follow the transfer operators by an orthogonal projection onto the new tangent space, 
\begin{equation}\label{eq:transfer_op_riemannian}
    \widetilde{\cI}_{h}^{H} = \P_{T_{X_H}\cMH} \circ \ \cI_{h}^{H}\big\lvert_{T_{X_{h}}\cMh}
    \quad \text{and} \quad
     \widetilde{\cI}_{H}^{h} = \P_{T_{X_h}\cMh} \circ \ \cI_{H}^{h}\big\lvert_{T_{X_{H}}\cMH}.
\end{equation}
This projection step is related to the so-called vector transport in retraction-based Riemannian optimization and can be seen as a first-order approximation of parallel transport in Riemannian geometry; see~\cite{AMS:2008}. As explained in the box below, the computation of the interpolation $\widetilde{\cI}_{h}^{H}$ exploits the factored form of tangent vectors. The implementation of the restriction $\widetilde{\cI}_{H}^{h}$ is similar and omitted.

\medskip
\fbox{\begin{minipage}{0.90\columnwidth}
\textbf{Interpolation} of $ \xi_{H} = U_{H} M_{H} V_{H} + U_{\mathrm{p},H} V_{H}\tr + U_{H} V_{\mathrm{p},H}\tr \in T_{X_H} \cMkH $\\
\textbf{Required}: $X_h = U_h \Sigma_h V_h\tr \in \cMkh$ and $X_{H} = U_H \Sigma_H V_H\tr \in \cMkH$
\begin{enumerate}[\qquad(1)]
\item Compute \textbf{factors}: \quad $ \widehat U_{\mathrm{p},h} = I_{H}^{h} U_{\mathrm{p},H}, \ \widehat M_{h} = M_{H}, \ \widehat V_{\mathrm{p},h} = I_{H}^{h} V_{\mathrm{p},H} $
\item \textbf{Normalize}: \quad $ U_{\mathrm{p},h} = \left( I - U_{h}U_{h}\tr\right) \widehat U_{\mathrm{p},h} $, $\ V_{\mathrm{p},h} = \left( I - V_{h}V_{h}\tr\right) \widehat V_{\mathrm{p},h} $ \\
\phantom{\textbf{Normalize}: \quad } $ M_{h} = U_{h}\tr \widehat U_{\mathrm{p},h} +  \widehat V_{\mathrm{p},h}\tr V_{h} + \widehat M_{h} $
\end{enumerate}
Result is $ \xi_{h} = U_{h} M_{h} V_{h} + U_{\mathrm{p},h} V_{h}\tr + U_{h} V_{\mathrm{p},h}\tr \in T_{X_h} \cMkh $ in the form~\eqref{eq:tan_vec_format_small_param}.
\end{minipage}}
\medskip

Like in \protect{\cite{Wen:2009}}, we will use injection and linear interpolation in the numerical experiments. In that case, the flop counts for computing $\cI_{h}^{H}$, $\widetilde{\cI}_{H}^{h}$, and $\widetilde{\cI}_{h}^{H}$ in factored form as explained above are both $O(nk^2 + k^3)$ for $\cMkh \subset \R^{n \times n}$.

\subsection{Smoothers}
In the context of optimization on manifolds, a smoother can be any cheap first-order optimization method for minimizing $ f_{h} $: given $ x_{h}^{(i)} $, it returns a tangent vector $ \xi_{h} $ such that, after retraction, the error of the new iterate $ \bar{x}_{h} = R_{x_{h}^{(i)}}(\xi_{h}) $ is smooth. In the Euclidean multilevel algorithm of~\cite{Wen:2009}, for example, a few steps of L-BFGS are used.

In our experiments, we simply use a fixed number of steps of Riemannian steepest descent; see~\cite{AMS:2008}.  In addition, we halve the step length found by the line-search method so that the resulting step better approximates one step of the Richardson iteration in linear multigrid.

\subsection{The Riemannian coarse-grid correction} 
Similar to Euclidean multilevel optimization, explained in~\cref{ref:sec_Euclidean_multilevel}, we also modify the objective function in the Riemannian setting.
To illustrate the generalization to the manifold case, let us first rewrite the Euclidean model \cref{eq:psi_H_CGC} as
\begin{equation}\label{eq:psi_cgc_Euclidean}
   \psi_{H}^\text{Euclidean}\colon \R^n \to \R, \quad  x_{H} \mapsto f_{H}(x_{H}) - g^{\mathrm{E}}(x_{H}, \kappa_{H} ),
\end{equation}
where $ x_{H} \coloneqq {x}_{H}^{(i)} + e_{H} $ is the full approximation. In the following, we describe how we turn this model into a function on manifolds. 

Let us assume that the algorithm at the coarse level starts at $ {x}_{H}^{(i)} \in \cMkH $. We consider as optimization variable a point on the manifold $ x_{H} \in \cMkH $. In the Riemannian setting, such a point $ x_{H} $ cannot be evaluated as in~\eqref{eq:psi_cgc_Euclidean} since the inner product $ g^{\mathrm{E}}(x_{H}, \kappa_{H}) $ is only defined for tangent vectors. We will therefore lift $ x_{H} $ to the tangent space at $ {x}_{H}^{(i)} $ by means of the inverse retraction when evaluating the inner product.\footnote{Recall from \cref{sec:ortho_retr} that this inverse is easy to compute for the orthographic retraction.} A coarse objective function suitable for Riemannian optimization is therefore given by
\begin{equation}\label{eq:CGC_on_cMk}
    \psi_{H} \colon \cMkH \to \R, \quad  x_{H} \mapsto f_{H}(x_{H}) -  g_{{x}_{H}^{(i)}}(R_{{x}_{H}^{(i)}}^{-1}(x_{H}),\kappa_{H}),
\end{equation}
where $ R_{{x}_{H}^{(i)}}^{-1} $ is the inverse retraction at $ {x}_{H}^{(i)} $, $ g_{{x}_{H}^{(i)}} $ denotes the Riemannian metric at $ {x}_H^{(i)} $, and $\kappa_{H} \in  T_{{x}_H^{(i)}}\cMkH $ is defined as
\begin{equation}\label{eq:kappa_H}
    \kappa_{H} = \grad f_{H} ( {x}_{H}^{(i)} ) - \widetilde{\cI}_{h}^{H}(\grad f_{h} (\bar{x}_{h})).
\end{equation}
Here, $\grad $ denotes the Riemannian gradient and $ \widetilde{\cI}_{h}^{H}(\grad f_{h} (\bar{x}_{h})) $ is the restricted Riemannian gradient coming from the fine-scale manifold. The restriction operator $ \widetilde{\cI}_{h}^{H} $ is defined as in \cref{eq:transfer_op_riemannian}, and the subtraction of the two tangent vectors is carried out in the factored format \cref{eq:tan_vec_format_small_param}. Let us denote by $ x_{H}^{(i+1)} $ the approximate minimizer of $ {\psi}_{H} $, and define the tangent vector $ {\eta}_{H} \coloneqq R_{{x}_{H}^{(i)}}^{-1}(x_{H}^{(i+1)}) $.

\subsection{Gradient of the coarse-grid model} During the optimization process, we need the Riemannian gradient of the coarse-grid correction function $\psi_H$. Recall from~\cref{eq:Riemannian_grad_as_projection} that this is simply the orthogonal projection of the Euclidean gradient onto the tangent space. 

To this end, let us simplify the notation by omitting $\cdot_H$ in~\eqref{eq:CGC_on_cMk}  to denote $\psi_H$ as
\begin{equation}\label{eq:CGC_on_cMk_2}
    \psi(x) = f(x) -  g_{x^{(i)}}(R_{x^{(i)}}^{-1}(x),\kappa),
\end{equation}
where $x$, $x^{(i)} \in \cMk $ and where the tangent vector $\kappa \in T_{ x^{(i)}} \cMk $  does not depend on $x$; see~\eqref{eq:kappa_H}. The only difficulty is thus the Euclidean gradient of the second term in~\eqref{eq:CGC_on_cMk_2}. Thanks to our choice of Riemannian metric on $ \cMk $, we have
\[
    g_{x^{(i)}}(R_{x^{(i)}}^{-1}(x),\kappa) = \langle \, R_{x^{(i)}}^{-1}(x), \, \kappa \, \rangle,
\]
where $ \langle \cdot, \cdot \rangle $ denotes the Frobenius inner product of two matrices.
By the chain rule, the Euclidean gradient of $ \langle \, R_{x^{(i)}}^{-1}(x), \, \kappa \, \rangle$ can therefore be written as the directional derivative
\[
    \nabla \langle \, R_{x^{(i)}}^{-1}(x), \, \kappa \, \rangle = \langle \, \nabla R_{x^{(i)}}^{-1}(x), \, \kappa \, \rangle = \D R_{x^{(i)}}^{-1}(x)[\kappa].
\]
For the orthographic retraction $R_x$, we know from~\cref{eq:inv_ortho_retr} that its inverse satisfies
\begin{equation*}
  R^{-1}_{x^{(i)}}(x) = \P_{T_{x^{(i)}}\cMk}(x) - x^{(i)}.
\end{equation*}
Since this is an affine linear function in $x$, its Fr\'{e}chet derivative is simply the orthogonal projection. We therefore obtain 
\[
    \D R_{x^{(i)}}^{-1}(x)[\kappa] = \P_{T_{x^{(i)}}\cMk}(\kappa) = \kappa,
\]
since $\kappa\in T_{x^{(i)}} \cMk $ by construction. 
Combining, we finally obtain the Riemannian gradient of $ {\psi} $ as
\[ 
    \grad{\psi}(x) = \P_{T_{x^{(i)}}\cMk}\!\big( \nabla f(x) \big) - \kappa.
\]

\begin{remark}
In the Euclidean multilevel optimization method from~\protect{\cite[eq.~(2.6)]{Wen:2009}}, an important property called \emph{first-order coherence} is introduced. In our Riemannian setting, it amounts to
\[
   g_{{x}_H} ( \grad{\psi}_H({x}_H), \xi_{H} ) = g_{{x}_h}  (\grad{f}_h(x_h), \xi_{h}),
\]
for any search direction $\xi_H \in T_{{x}_H} \cMH$ with $x_H =  \cI_{h}^{H}(x_h)$ and $\xi_h =   \widetilde{\cI}_{H}^{h} (\xi_H)$. This is a desirable property since it ensures the same slope of the objective functions on the fine and coarse grids. Practically, this equation imposes a relation between the intergrid transfer operators in the multilevel algorithm. In our setting as explained in~\cref{sec:transfer_operators}, one can show that it requires $ I_{H}^{h} = (I_{h}^{H})\tr $. This is indeed a typical choice in multigrid algorithms. It is, for example, satisfied for the injection  $I_{H}^{h}$ and linear interpolation $I_{h}^{H}$.
\end{remark}

\subsection{Final algorithm: Riemannian multigrid line search}\label{sec:RMGLS_code} 
In the following box, we have listed the final Riemannian multigrid line-search algorithm to optimize an objective function on a Riemannian manifold. The smoother is denoted by the function $ \mathrm{SMOOTH}$ and corresponds to $\nu_1$ or $\nu_2$ steps of steepest descent for $f_h$.

\medskip

\begin{center}
\noindent
\fbox{\begin{minipage}{0.95\columnwidth}
\textbf{One RMGLS iteration} starting at $x_{h}^{(i)}$ to minimize $f_h$.

\begin{enumerate}[\qquad(1)]
\item \textbf{Pre-smoothing}: $ \bar{x}_{h} = \mathrm{SMOOTH}^{\nu_{1}}( x_{h}^{(i)}, f_{h} ) $

\item \textbf{Coarse-grid correction}: 

\begin{enumerate}[(a)]

\item \textbf{Restrict} as in section~\ref{sec:transfer_operators} to the coarse manifold: $ x_{H}^{(i)} = \cI_{h}^{H} ( \bar{x}_{h} ) $

\item Compute the \textbf{linear correction term}: 
\[
 \kappa_{H} = \grad f_{H} ( {x}_{H}^{(i)} ) - \widetilde{\cI}_{h}^{H}(\grad f_{h} (\bar{x}_{h})) 
\]   

\item Define the \textbf{coarse-grid objective function} 
\[
\psi_H (x_{H}) = f_{H}(x_{H}) -  g_{{x}_{H}^{(i)}}(R_{{x}_{H}^{(i)}}^{-1}(x_{H}),\kappa_{H})
\]

\item Compute an \textbf{approximate minimizer $x_{H}^{(i+1)}$} starting at $x_{H}^{(i)}$ to minimize $\psi_H$ using either
\begin{itemize}
\item a Riemannian trust-region method (if $\cMH$ is small)
\item one recursive RMGLS iteration (otherwise)
\end{itemize}

\item Compute the \textbf{coarse-grid correction}: $\eta_{H} = R_{x_{H}^{(i)}}^{-1} (x_{H}^{(i+1)})$

\item \textbf{Interpolate} as in section~\ref{sec:transfer_operators} to the fine manifold: $ \eta_{h} = \widetilde{\cI}_{H}^{h}(\eta_{H}) $

\item Compute the \textbf{corrected approximation} on the fine manifold:
 \[
   \widehat{x}_{h} = R_{\bar{x}_{h}}(\alpha^{\ast} {\eta}_{h}) \quad \text{with $\alpha^{\ast}$ obtained from  line search}
\]

\end{enumerate}

\item \textbf{Post-smoothing}: $ x_{h}^{(i+1)} = \mathrm{SMOOTH}^{\nu_{2}}(\widehat{x}_{h}, f_{h} ) $

\end{enumerate}

\end{minipage}}
\end{center}

\medbreak
\begin{remark}
The RMGLS algorithm above is very similar to the way one FAS multigrid iteration is presented in~\cite[p.~157]{Trottenberg:2000}.
\end{remark}

\begin{remark} For efficiency reasons, it is crucial to implement the algorithm without forming full matrices, i.e., always exploiting the low-rank format explicitly, also when evaluating the objective function $f$. More details will be given in~\cref{{sec:variational_problems}}.
\end{remark}

\section{A more accurate line search}
\label{sec:linesearch}

Our optimization algorithm RMGLS performs line searches during the smoothing steps and the application of the coarse-grid correction. Like other optimization methods that only use first-order information, the convergence of RMGLS to the stationary point is generically linear. We will explain below that this makes it difficult to achieve high accuracy in finite precision arithmetic when using standard line searches, like the weak Wolfe conditions. A more accurate line search was proposed by~\cite{Hager:2005} in the context of a new nonlinear CG method. Here we explain how to adapt this line-search method to the Riemannian setting.

\subsection{Inaccuracy in standard line search}

\begin{figure}[htbp]
   \centering
   \includegraphics[width=0.50\columnwidth]{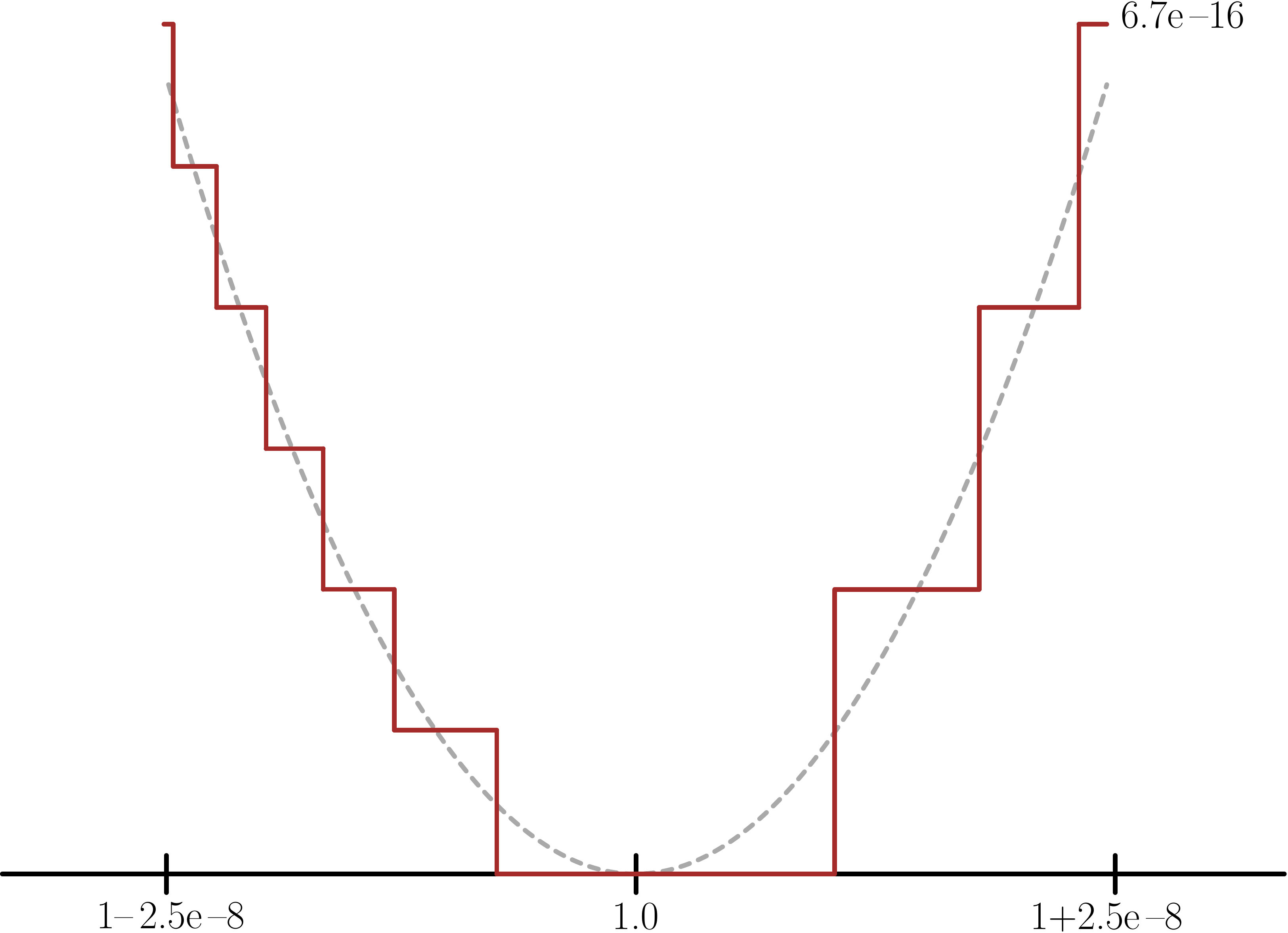}
   \caption{Exact and numerical graphs of $ f(x) = 1-2x + x^{2} $ near $ x=1 $ (adapted from \protect{\cite[§4]{Hager:2005}}). The dotted line is the exact $f$, while the solid line is its representation in double precision with $ \epsmach \approx 10^{-16}$.}\label{fig:quadratic_numerical_exact}
\end{figure}

The usual stopping criterion for line search is the weak Wolfe conditions, which we recall here.
Let $ f $ be a differentiable objective function. Let $ x $ be the current iterate, $ g = \nabla f(x) $ the gradient, and $ d $ the search direction.  
The weak Wolfe conditions for the step size $ \alpha > 0$ are defined by
\[
   f(x + \alpha \cdot d ) - f(x) \leq \delta \alpha \, d \tr \nabla f(x), \qquad  d\tr \nabla f(x+ \alpha \cdot d) \geq \sigma \, d \tr \nabla f(x),  
\]
with  $0 < \delta \leq \sigma < 1$. The first inequality is known as sufficient decrease, or Armijo, condition, while the second represents a curvature condition.
Observe that the weak Wolfe conditions can also be recast in terms of $ \phi(\alpha) \coloneqq f(x + \alpha\cdot d)$ as follows:
\begin{equation}\label{eq:weak_wolfe}
   \delta\, \phi'(0) \geq \frac{\phi(\alpha) - \phi(0)}{\alpha},  \qquad  \phi'(\alpha) \geq \sigma \, \phi'(0), 
\end{equation}
with $0 < \delta \leq \sigma < 1$. In finite precision, the weak Wolfe conditions can be difficult to satisfy very accurately due to roundoff error when $x$ is very close to the local minimum of $f$. For a smooth objective function with a strict local minimum, the function $f$ is locally quadratic and its minimum can indeed only be determined within $ \sqrt{\epsmach} $, with $ \epsmach $ the machine epsilon; see \cref{fig:quadratic_numerical_exact}. 

\subsection{Approximate Wolfe conditions}
To prevent the loss of accuracy during standard line search, Hager and Zhang proposed in \protect{\cite[§4]{Hager:2005}} to relax the weak Wolfe conditions~\eqref{eq:weak_wolfe} to the so-called \textit{approximate} Wolfe conditions.
Their main observation is that, in a neighborhood of a local minimum, the first condition in \cref{eq:weak_wolfe} is difficult to satisfy since $\phi(\alpha) \approx \phi(0)$ which makes the subtraction $ \phi(\alpha) - \phi(0) $ relatively inaccurate \protect{\cite[§3]{Hager:2006}}.
This leads them to introduce the approximate Wolfe conditions \protect{\cite[eq.~(4.1)]{Hager:2005}}
\begin{equation}\label{eq:approximate_weak_wolfe}
   \left( 2\delta - 1\right) \phi'(0) \geq \phi'(\alpha) \geq \sigma \, \phi'(0), \qquad 0 < \delta < 0.5, \quad \delta \leq \sigma < 1.
\end{equation}
Here, the first inequality is an approximation of the first condition in \cref{eq:weak_wolfe}, but the second inequality coincides with the second condition in \cref{eq:weak_wolfe}.\footnote{It would therefore be more appropriate to talk about the \textit{approximate Armijo} condition than the \textit{approximate Wolfe} conditions, but we stick with the latter name as in \cite{Hager:2005,Hager:2006}.}

The idea behind this approximation comes from replacing $\phi$ by its quadratic interpolant $q$ that satisfies $ q(0) = \phi(0) $, $ q'(0) = \phi'(0) $, and $ q'(\alpha) = \phi'(\alpha) $. Now the finite difference quotient in the first Wolfe condition can be approximated as
\[
   \frac{\phi(\alpha) - \phi(0)}{\alpha} \approx \frac{q(\alpha) - q(0)}{\alpha} = \frac{\phi'(\alpha)+ \phi'(0)}{2}. 
\]
With this approximation, the subtraction $ q(\alpha) - q(0) $ can be computed more accurately as $\phi'(\alpha)+ \phi'(0)$, thereby circumventing the possible cancellation due to roundoff errors in the original difference $ \phi(\alpha) - \phi(0) $.

\subsection{The Hager--Zhang bracketing}
\label{sec:hzls}
Another component in the Hager--Zhang line-search method is a new procedure to determine a step length $\alpha$ that satisfies the approximate Wolfe conditions~\eqref{eq:approximate_weak_wolfe}. It combines the secant and bisection methods, as in other line-search methods, but applied to find a zero of the derivative $\phi'$ instead of a minimum of $\phi$. In particular, the method from \protect{\cite[§3]{Hager:2006}} will generate a nested sequence of bracketing intervals that are guaranteed to contain an acceptable step length $\alpha$. A typical interval $[a,b]$ in this sequence satisfies the opposite slope condition
\[
   \phi'(a) < 0, \qquad \phi'(b) \geq 0.
\]
This is nothing else than the opposite sign condition of the bisection method translated to the derivative, meaning that the derivative changes sign in the bracketing interval.

Since finding a zero of an almost linear function is better conditioned numerically than finding a minimum of an almost quadratic, the Hager--Zhang bracketing search is more accurate. In addition, it also explicitly takes into account roundoff error when $\phi'(\alpha) \approx 0$. We omit the many technical details and refer to \protect{\cite[§3]{Hager:2006}}.

\subsection{Numerical example}
\begin{figure}[htbp]
  \centering
  \includegraphics[width=0.65\textwidth]{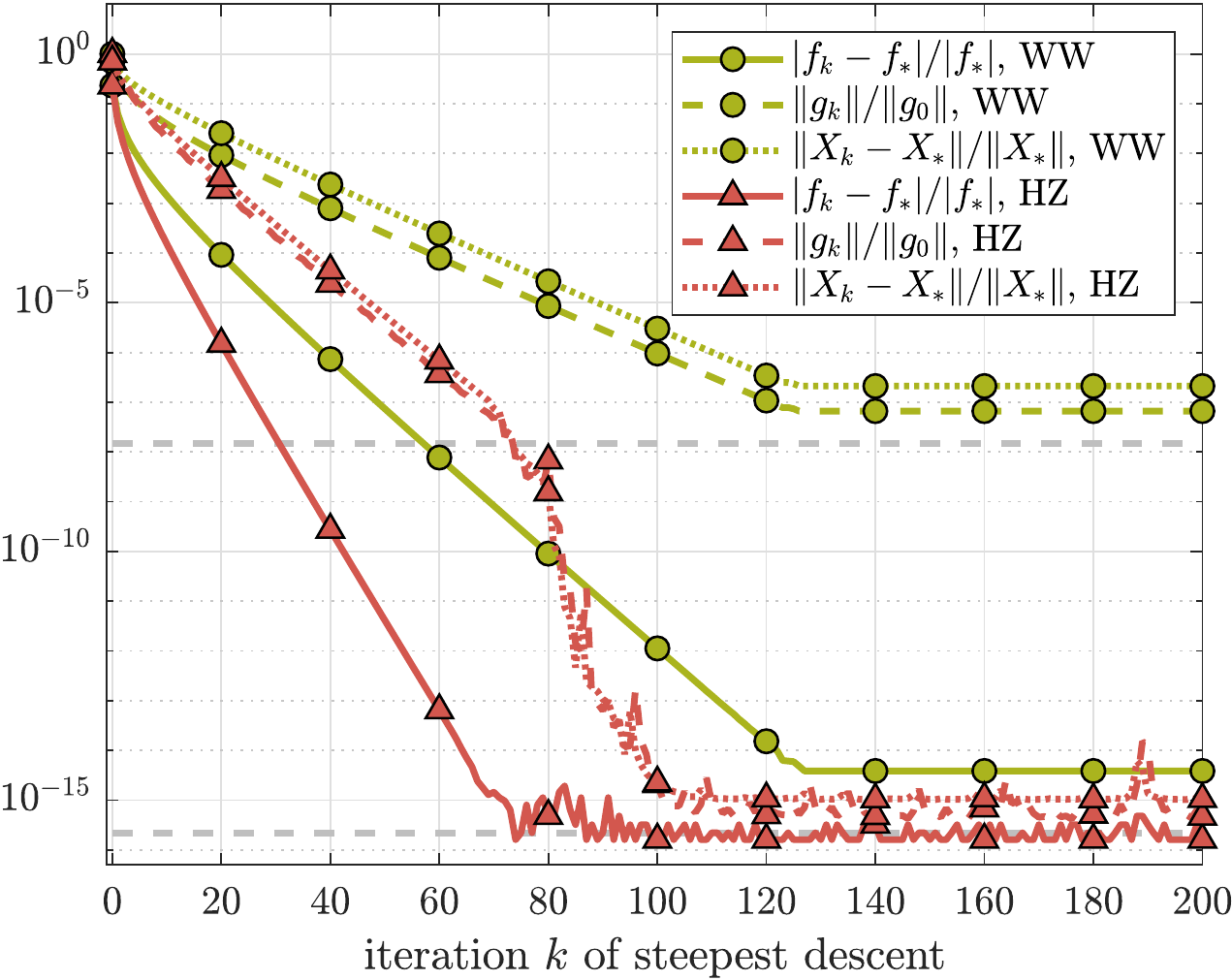}
  \caption{Convergence behavior of line search with weak Wolfe (WW) or Hager--Zhang (HZ) when applied to a quadratic function $f$. The objective function is denoted by $f_k$ and the gradient by $g_k$. The horizontal dashed lines indicate $\sqrt{\epsmach}$ and $\epsmach$.}
  \label{fig:quadratic_matrix_hz_vs_ww}
\end{figure}

To illustrate the convergence behavior of the Hager--Zhang line search, consider the following numerical example in which steepest descent is used to minimize the quadratic cost function $ f \colon \R^{n \times n} \to \R $, defined by
\[
   f(X) = \tfrac{1}{2} \trace(X\tr A X) - \trace( X\tr B ).
\]
We take $ n = 100 $, the condition number of the symmetric positive definite matrix $A$ as $ \kappa(A) = 10 $, and $ B = A X^{\ast} $, where $ X^{\ast} $ is the exact solution to the problem $ A X = B $. The starting point of the optimization is a random initial guess $X^{(0)}$. We compare the results of steepest descent using weak Wolfe conditions and the Hager--Zhang line search.
From \cref{fig:quadratic_matrix_hz_vs_ww}, we see that the gradient norm stagnates at about $10^{-8}$ for the weak Wolfe conditions. In contrast, the approximate Wolfe conditions used by the Hager--Zhang line search allow to reach an accuracy on the order of $ \epsmach $ ($ \approx  10^{-16} $ in double precision) in both the objective value and the gradient norm\footnote{Another line search variant can be obtained by introducing the approximate Armijo condition in \textsc{Matlab}'s native \texttt{fzero} to find directly a zero of the derivative $\phi'$. From the numerical experiments it appears that this version is much less efficient than the Hager--Zhang line search, even though it attains the same accuracy level.}. 
Observe that the error for the objective function is still of full accuracy because it is equivalent to the square of $ \| X_{k} - X_{\ast} \| $. On the other hand, to obtain a small value for $ \| X_{k} - X_{\ast} \|/\|X_{\ast}\|$, a small objective value alone is not sufficient and a small gradient is also needed.

Let us also point out that, for this example, when using the approximate Wolfe conditions the number of function evaluations is about 55\% less than the one attained by using the weak Wolfe conditions. This is likely because the standard line search wastes a lot of effort in bracketing the function $\phi(\alpha)$ that becomes noisy due to roundoff error when $\alpha$ is close to a stationary point.

\subsection{Riemannian Hager--Zhang line search}
\label{sec:riemannian_hzls}

The Hager--Zhang line search explained above can be readily extended to Riemannian manifolds by applying it to the retracted objective function $\phi(\alpha) = f(R_x(\alpha \cdot \eta))$ along the search direction  $\eta \in T_x \cM$   with $\alpha \geq 0$ the step length. We call this generalization Riemannian Hager--Zhang line search.  The only difficulty in its implementation is the need for $\phi'$ since this requires computing the derivative $\mathrm{d} R_x(\alpha \cdot \eta) / \mathrm{d}t$, which is cumbersome for general retractions $R_x$. Observe that line searches based on the Armijo condition alone only require this derivative at $\alpha=0$ for which it equals $\eta$ by definition of a retraction.

Fortunately, in Riemannian optimization we can choose a retraction that better suits our needs. As mentioned above, we have chosen the orthographic retraction for the manifold of fixed-rank matrices $ \cMk $ because it has an explicit inverse; see \cref{sec:ortho_retr}.
Let us show that its derivative can also be efficiently calculated.

Let $X \in \cMk$ and $R_X$ the orthographic retraction. Recall that $f \colon \cMk \to \R$. By the chain rule, we get for $\phi(t) = f(R_X(t\eta))$ that
\begin{equation}\label{eq:der_phi_t}
   \phi'(t) = \big\langle \nabla f(R_{X}(t\eta)), \tfrac{\mathrm{d} }{\mathrm{d} t} R_{X}(t\eta)\big\rangle = \trace\!\Big( \nabla f(R_{X}(t\eta)) \tr \tfrac{\mathrm{d}}{\mathrm{d} t} R_{X}(t\eta) \Big),
\end{equation}
where $ \nabla f $ is the Euclidean gradient of $ f $. Using~\eqref{eq:ortho_graphic_R}, we can work out the standard derivative
\begin{align*}
   \begin{split}
      \tfrac{\mathrm{d} }{\mathrm{d} t}R_{X}(t\eta) &= \big( U (\Sigma + tM) + t U_{\mathrm{p}} \big) (\Sigma+ t M)^{-1} (MV\tr + V_{\mathrm{p}}\tr) \\
          & - \big( U(\Sigma + t M) + t U_{\mathrm{p}} \big) (\Sigma + t M)^{-1} M (\Sigma+tM)^{-1} \big( (\Sigma+t M)V\tr + t V_{\mathrm{p}}\tr \big) \\
          & + ( U_{\mathrm{p}} + UM ) (\Sigma + t M)^{-1} \big( (\Sigma+t M)V\tr + t V_{\mathrm{p}}\tr \big),
   \end{split}
\end{align*}
where $ X = U\Sigma V\tr $ and $  \eta = UMV\tr + U_{\mathrm{p}}V\tr + UV_{\mathrm{p}}\tr $ as in \cref{eq:tan_vec_format_small_param}.

In \cref{eq:der_phi_t}, we need to evaluate $ \trace(A\tr\! B) $. For computational efficiency, we want to avoid the naive multiplication $ A\tr\! B $ since it costs $ O(n^{3}) $ flops. A more efficient approach is to rewrite  the derivative in the factorized format $ \tfrac{\mathrm{d} }{\mathrm{d} t}R_{X}(t\eta) = GH\tr $ by defining
\[
   G = \left[\begin{array}{ccc}
          -( U + t U_{\mathrm{p}} (\Sigma + t M)^{-1}) M (\Sigma+tM)^{-1}  & \ \ U + t U_{\mathrm{p}} (\Sigma + t M)^{-1}  & \ \ U_{\mathrm{p}} + UM
       \end{array}\right]
\]
and
\[
   H = \left[\begin{array}{ccc}
          V(\Sigma + t M)\tr + tV_{\mathrm{p}} & \quad VM\tr + V_{\mathrm{p}} & \quad V + tV_{\mathrm{p}}(\Sigma+tM)^{-\textsf {T}}
       \end{array}\right].
\]
Observe that $G,H \in \R^{n \times 3k}$. Assuming a similar factorization for $ \nabla f(R_{X}(t\eta)) = \widetilde{G} \widetilde{H}\tr $ with $\widetilde{G}, \widetilde{H} \in \R^{n \times \widetilde{k}}$, the trace in~\eqref{eq:der_phi_t} can then be computed as
\[
   \phi'(t) = \trace\!\big( \widetilde{H} \widetilde{G}\tr\! GH\tr \big) = \trace\!\big( (\widetilde{G}\tr\! G) (H\tr\!\widetilde{H}) \big)
\]
at a cost of $O((n+k)k\widetilde{k} )$. In typical applications targeting low-rank approximations, $\widetilde{k}$ is larger than $k$ but significantly smaller than $n$. For example, in our numerical experiments below, $\widetilde{k} = O(k^2)$ showing a large reduction from $O(n^3)$ when $k$ is small.

\section{Numerical experiments for two variational problems}
\label{sec:variational_problems}
We report on numerical properties of the proposed algorithm, RMGLS, by applying it to the variational problems presented in this section. These are large-scale finite-dimensional optimization problems arising from the discretization of infinite-dimensional problems. Because of their underlying PDEs, these variational problems present a natural multilevel structure. Variational problems of this type have been considered as benchmarks in other nonlinear multilevel algorithms \protect{\cite{Henson:2003,Gratton:2008,Wen:2009}}. For the theoretical aspects of variational problems, some good references are \protect{\cite{Brenner:2007,LeDret:2016}}.

The experiments below were performed by recursively executing RMGLS in a V-cycle manner for both problems, as explained in~\cref{sec:RMGLS_code}. Unless otherwise noted, the Riemannian version of the Hager--Zhang line search was used. The algorithm was implemented in \textsc{Matlab} and is publicly available \protect{\cite{Sutti:2020}}.

\subsection{A linear problem (Lyapunov equation)}
\label{sec:var_pb_1}
We consider the minimization problem
\begin{equation}\label{eq:var_pb_1}
   \begin{cases}
      \displaystyle\min_{w} \cF(w(x,y)) = \int_{\Omega} \tfrac{1}{2} \| \nabla w(x,y) \|^{2} - \gamma(x,y)\,w(x,y) \dx\dy \\
      \quad \text{such that} \quad w=0 \ \text{on} \ \partial\Omega,
   \end{cases}
\end{equation}
where $ \nabla = \big( \frac{\partial}{\partial x}, \frac{\partial}{\partial y} \big) $, $ \Omega = [0,1]^{2} $ and $ \gamma $ is the source term.
The variational derivative (Euclidean gradient) of $ \cF $ is
\begin{equation}\label{eq:var_pb_1_grad}
   \frac{\delta \cF}{\delta w} = -\Delta w - \gamma.
\end{equation}
A critical point of~\eqref{eq:var_pb_1} is thus also a solution of the elliptic PDE $ -\Delta w = \gamma$.

\subsubsection{Discretization of the objective function}
\label{sec:var_pb_1_discr_functional}
We use a standard finite difference discretization for~\eqref{eq:var_pb_1}. In particular, 
$\Omega$ is represented at level $\ell$ as a square grid
\begin{equation*}
   \Omega_{\ell} = \lbrace (x_{i},y_{j}) \mid x_{i} = ih_{\ell}, \ y_{j} = jh_{\ell}, \ i = 0,1, \ldots, n_{\ell}, \ j = 0,1,\ldots,n_{\ell} \rbrace, \quad n_{\ell} = 2^{\ell}, 
\end{equation*}
yielding a square mesh of uniform mesh width $ h_{\ell} = 1/n_{\ell}$. The unknown $w$ on $\Omega_\ell$ is denoted by  $ w_{ij} \coloneqq w(x_{i}, y_{j}) $, and likewise for $ \gamma_{ij} \coloneqq \gamma(x_{i}, y_{j}) $, where we have omitted the dependence on $\ell$ in the notation for readability. 
The partial derivatives are discretized as forward finite differences
\begin{equation}\label{eq:forward_difference}
   \partial w_{x_{ij}} = \tfrac{1}{h} (w_{i+1,j} - w_{i,j}), \qquad \partial w_{y_{ij}} = \tfrac{1}{h}(w_{i,j+1} - w_{i,j}).
\end{equation}
The discretized version of $\cF$ therefore becomes
\begin{equation}\label{eq:var_pb_1_discr_functional}
   \cF_{h} = h^{2} \sum_{i,j=0}^{2^{\ell} - 1} \big( \tfrac{1}{2} ( \partial w_{x_{ij}}^{2} + \partial w_{y_{ij}}^{2} ) - \gamma_{ij} \, w_{ij} \big).
\end{equation}

In order to find a low-rank approximation of~\eqref{eq:var_pb_1} with RMGLS, the unknown $w_{ij}$ from above will be approximated as the $ij$th entry of a matrix $W_{h} \in \R^{n \times n}$ of rank $k$. For efficiency, this matrix is always represented in the factored form $ W_{h} = U\Sigma V\tr $. Likewise, we gather all $\gamma_{ij}$ in a factored matrix $\Gamma_h= U_\gamma \Sigma_\gamma V_\gamma \tr$ of rank $k_\gamma$. In the experiments below, $k_\gamma=5$.

For reasons of computational efficiency, it is important to exploit these low-rank forms in the execution of RMGLS. For the objective value $ \cF_{h} $ this can be done as follows. First, observe that the first term satisfies
\begin{equation}\label{eq:first_term}
   I \coloneqq \sum_{i,j=0}^{2^{\ell} - 1}  ( \partial w_{x_{ij}}^{2} + \partial w_{y_{ij}}^{2} ) = \| \partial W_{x} \|^{2}_{\F} + \| \partial W_{y} \|^{2}_{\F},
\end{equation}
where $\partial W_{x}$, $ \partial W_{y} \in \R^{n \times n}$ contain the derivatives $ \partial w_{x_{ij}} $, $ \partial w_{y_{ij}} $.
Then it is easy to verify from~\eqref{eq:forward_difference} that 
\begin{equation*}
   \partial W_{x} = L W_{h}  \quad \text{and} \quad  \partial W_{y} = W_{h} L\tr \quad \text{with} \quad
   L = \frac{1}{h}
   \begin{bmatrix}
       -1 &   1 &          &         &    \\
          &  -1 &        1 &         &    \\ 
          &     &   \ddots &  \ddots &    \\
          &     &          &      -1 &  1
   \end{bmatrix}.
\end{equation*}
Substituting this factorization and $ W_{h} = U \Sigma V\tr $ in~\cref{eq:first_term}, we get
\begin{equation*}
I = \| (LU) \Sigma \|^{2}_{\F} + \| (LV) \Sigma \|^{2}_{\F}.
\end{equation*}
To recast the second term in~\cref{eq:var_pb_1_discr_functional} using matrices, observe that
\[
   I\! I \coloneqq \sum_{i,j} \gamma_{ij} \, w_{ij} = \sum_{i,j}(\Gamma_{h} \odot W_{h})_{ij} = \trace( \Gamma_{h}\tr W_{h} ) = \trace\!\big( \Sigma_{\gamma} (U_{\gamma}\tr U)\Sigma (V\tr V_{\gamma}) \big),
\]
where $ \odot $ denotes the elementwise (or Hadamard) product of two matrices. Summing the terms $I$ and $I\! I$, we finally obtain
\[
    \cF_{h} = \tfrac{h^{2}}{2} \Big( \| (LU) \Sigma \|^{2}_{\F} + \| (LV) \Sigma \|^{2}_{\F} -2 \trace\!\big( \Sigma_{\gamma} (U_{\gamma}\tr U)\Sigma (V\tr V_{\gamma}) \big) \Big),
\]
which can be evaluated in $O( (n+k_{\gamma})k_\gamma k)$ flops.

\subsubsection{Discretization of the gradient}
\label{sec:var_pb_1_grad_discretized}
The discretization of \cref{eq:var_pb_1_grad} gives
\begin{equation}\label{eq:var_pb_1_grad_discretized}
    G_{h} = h^{2} \left( A W_{h} + W_{h} A - \Gamma_{h} \right),
\end{equation}
where $ A $ is the discretization of $ -\Delta $ by a second-order central difference, i.e.,
\begin{equation}\label{eq:var_pb_1_Laplacian_A}
   A = \frac{1}{h^{2}}
   \begin{bmatrix}
    2 &      -1 &         &         &     \\
   -1 &       2 &     -1  &         &     \\ 
      &  \ddots & \ddots  &  \ddots &     \\
      &         &     -1  &       2 &  -1 \\
      &         &         &      -1 &   2
   \end{bmatrix}.
\end{equation}
Observe that $G_h=0$ in \cref{eq:var_pb_1_grad_discretized} --- and hence for $W_{h}$ a critical point of~\eqref{eq:var_pb_1_discr_functional} --- coincides with a solution to the Lyapunov equation $A W_{h} + W_{h} A = \Gamma_{h}$.

Like for the discretized objective function above, we represent the discretized gradient $G_h$ as a factored matrix. Using the same notation as above, this can be done as follows:
\begin{align*}
   G_{h} & = h^{2} \left( AU\Sigma V\tr + U\Sigma V\tr\!\!A - \Gamma_{h} \right) \\
   & = h^{2} \left( (AU)\Sigma V\tr + U \Sigma(V\tr\!\!A) - U_{\gamma} \Sigma_{\gamma} V_{\gamma}\tr \right) \\
   & = h^{2}
   \begin{bmatrix}
      AU  &  U  &  U_{\gamma}
   \end{bmatrix}
   \blkdiag\!\left(\Sigma, \Sigma, -\Sigma_{\gamma} \right)
   \begin{bmatrix}
      V  &  AV  &  V_{\gamma} 
   \end{bmatrix}\tr,
\end{align*}
where $ \blkdiag\!\left(\Sigma, \Sigma, -\Sigma_{\gamma} \right) $ is the block diagonal matrix created by aligning the matrices $\Sigma$, $\Sigma$, and $-\Sigma_{\gamma}$ along the main diagonal. The gradient $G_h$ can be represented in only $O(nk)$ flops for computing $AU$ and $AV$.

We introduce the notation $ \xi_{h} $ for the Riemannian gradient and recall that it is given by the projection~\eqref{eq:Riemannian_grad_as_projection}
\[
   \xi_{h} = \P_{T_{W_{h}}\cMkh} (G_{h}).
\]
Its norm $ \| \cdot \|_{\F} $ can be directly computed from the format \cref{eq:tan_vec_format_small_param} as
\[
   \| \xi_h \|_{\F} = \sqrt{\| M \|_{\F}^{2} + \| U_{\mathrm{p}} \|_{\F}^{2} + \| V_{\mathrm{p}} \|_{\F}^{2}}.
\]

\subsubsection{Numerical results}
\label{sec:num_res_var_pb_1}

As mentioned above, the unconstrained minimizer of $\cF_h$ over $\R^{n \times n}$ is also a solution of a Lyapunov equation. Restricted to $\cMkh$ and for small $k_\gamma$, this is a typical benchmark problem for low-rank methods; see, e.g., \protect{\cite[§4.4]{Simoncini:2016}}. In particular, it guarantees the existence of an approximation of rank
\[
   k = O\big(\log(1/\varepsilon) \, \log(\kappa(A)) \, k_\gamma \big),
\]
with error at most  $ \varepsilon $ and $ \kappa(A) $ the condition number of $A$. In the experiments, we have $k_\gamma = 5$ and
\begin{equation}\label{eq:rank5_rhs}
   \gamma(x,y) = e^{x-2y} \, \sum_{j=1}^{5} 2^{j-1} \sin( j\pi x) \sin(j \pi y).
\end{equation}

We now report on the behavior of RMGLS. In all cases, we used 5 pre- and 5 post-smoothing steps, and coarsest scale $\ell_{\mathrm{c}} = 5$. To monitor the convergence behavior of RMGLS, we have considered three quantities. In all formulas, ${\cdot}^{(i)}$ indicates that a quantity was evaluated at the $i$th outer iteration of RMGLS.
\begin{itemize}
\item[(a)] The relative error of the discretized objective function $\cF_h$:
\[
   \textrm{err-$\cF$}(i) \coloneqq  \lvert \cF^{(i)}_{h} - \cF^{(\ast)}_{h} \rvert/\lvert \cF^{(\ast)}_{h} \rvert.
\]
Here, $ \cF^{(\ast)}_{h} $ is the minimal value over $\cMkh$ of the original objective function in~\eqref{eq:var_pb_1}. It is approximated by minimizing $ \cF_{h}$ on $\cMkh$ with the Riemannian trust-region (RTR) method~\cite{AMS:2008}, terminated when the Riemannian gradient norm is smaller than $ 10^{-13}$.
\item[(b)] The Frobenius norm of the normalized Riemannian gradient:
\[
   \textrm{R-grad}(i) \coloneqq \| \xi_{h}^{(i)} \|_{\F}/\| \xi_{h}^{(0)} \|_{\F}.
\]
\item[(c)] The relative error in Frobenius norm of the low-rank approximation:
\[
   \textrm{err-$W$}(i) \coloneqq \| W_{h}^{(i)} - W^{(\ast)}_{h} \|_{\F}/\| W^{(\ast)}_{h} \|_{\F}.
\]
Here, $ W^{(\ast)}_{h} $ is the minimizer of $\cF_h$ over $\R^{n \times n}$. It is computed with a Euclidean trust-region method, terminated when the Euclidean gradient norm is smaller than $ 10^{-14}$. No rank truncation was used for $ W^{(\ast)}_{h} $ and no problem with multiple local minima occurred.\footnote{For the linear problem, we can of course also directly solve the Lyapunov equation. However, this is not feasible for the nonlinear problem below.}
\end{itemize}

\begin{figure}[htbp]
  \centering
  \includegraphics[width=0.65\textwidth]{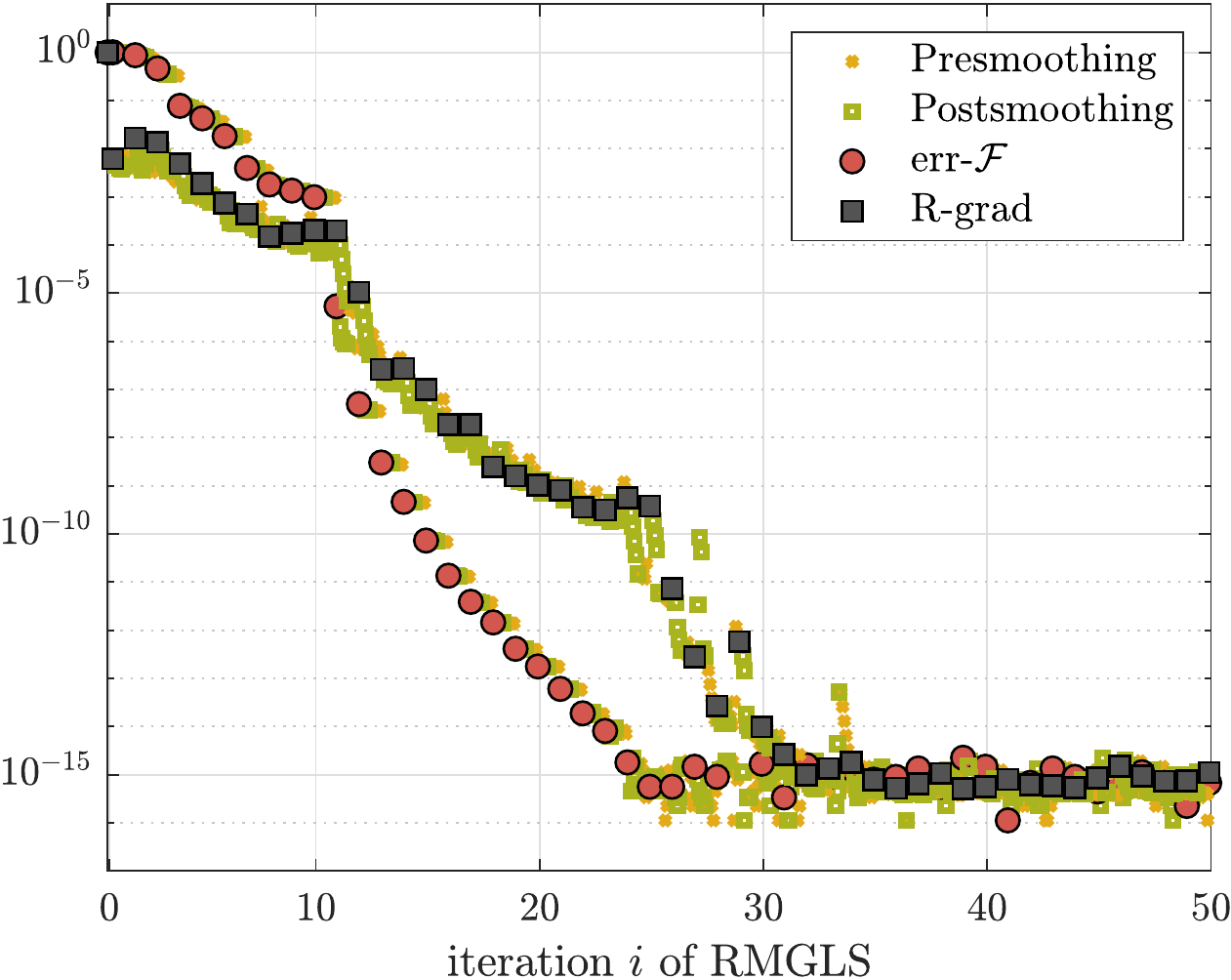}
  \caption{Convergence of $ \mathrm{err}\textrm{-}\cF $ and $ \mathrm{R\textrm{-}grad} $ for level $ \ell_{\mathrm{f}} = 8 $ and rank $ k = 5 $, for the problem of~\cref{sec:var_pb_1}.}\label{fig:LYAP_CostGrad_5_8_K5_HZ_SS5_mgiter_50}
  \end{figure}

In~\cref{fig:LYAP_CostGrad_5_8_K5_HZ_SS5_mgiter_50}, the convergence of the objective function and gradient norm are depicted for RMGLS with finest scale $ \ell_{\mathrm{f}} = 8 $ and rank $ k = 5 $. We observe that the objective function has converged already after 25 iterations, whereas the gradient norm continues to decrease until iteration 35. This difference indicates that using a stopping criterion based on the objective function alone can be misleading if we want the most accurate stationary point, and it is better to use a criterion based on the gradient norm.

\begin{figure}[htbp]
    \centering
    \begin{minipage}{0.48\textwidth}
        \centering
        \includegraphics[width=\textwidth]{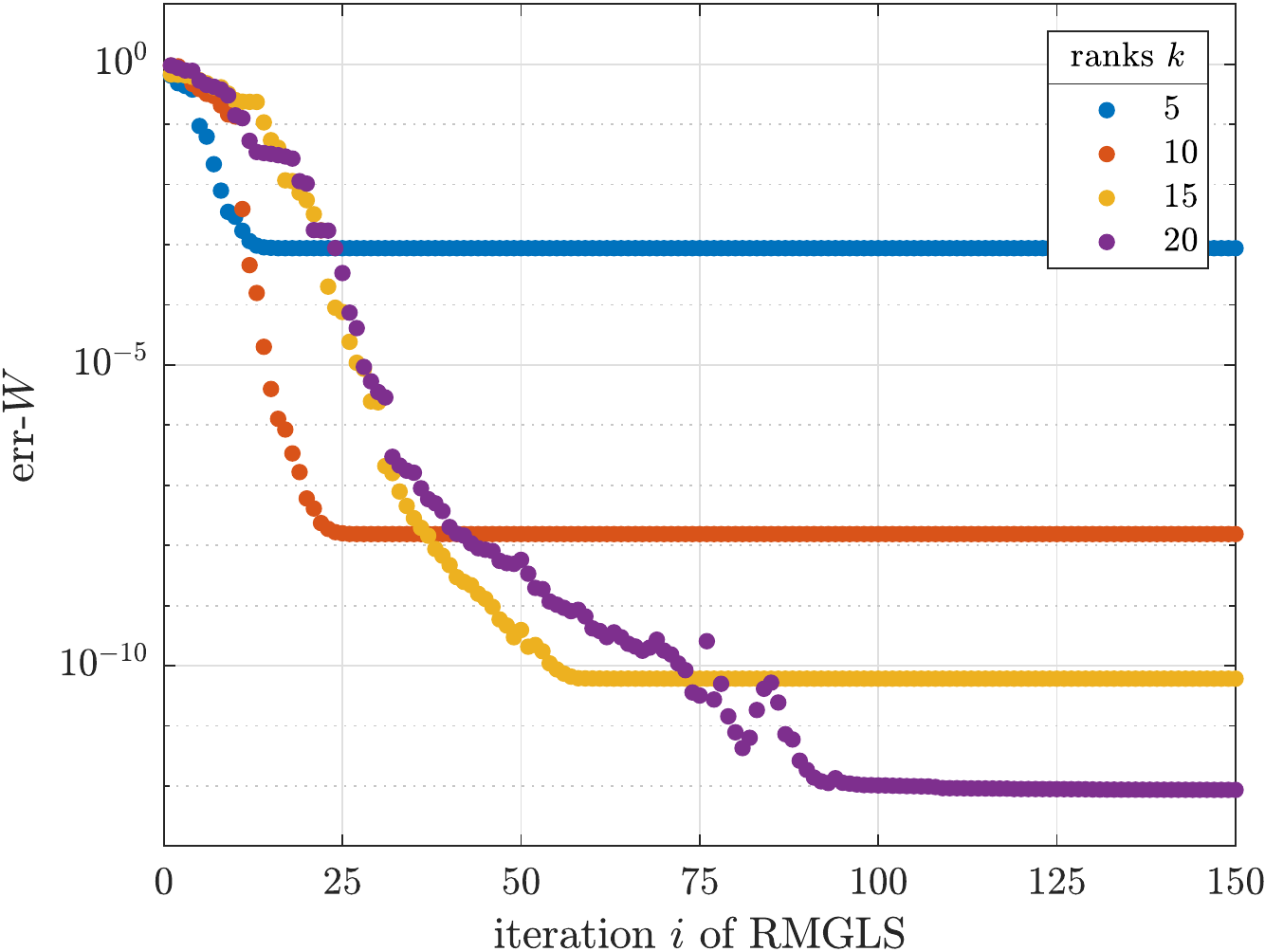}
        {\scriptsize (a) \emph{Hager--Zhang}}
    \end{minipage}\hfill
    \begin{minipage}{0.48\textwidth}
        \centering
        \includegraphics[width=\textwidth]{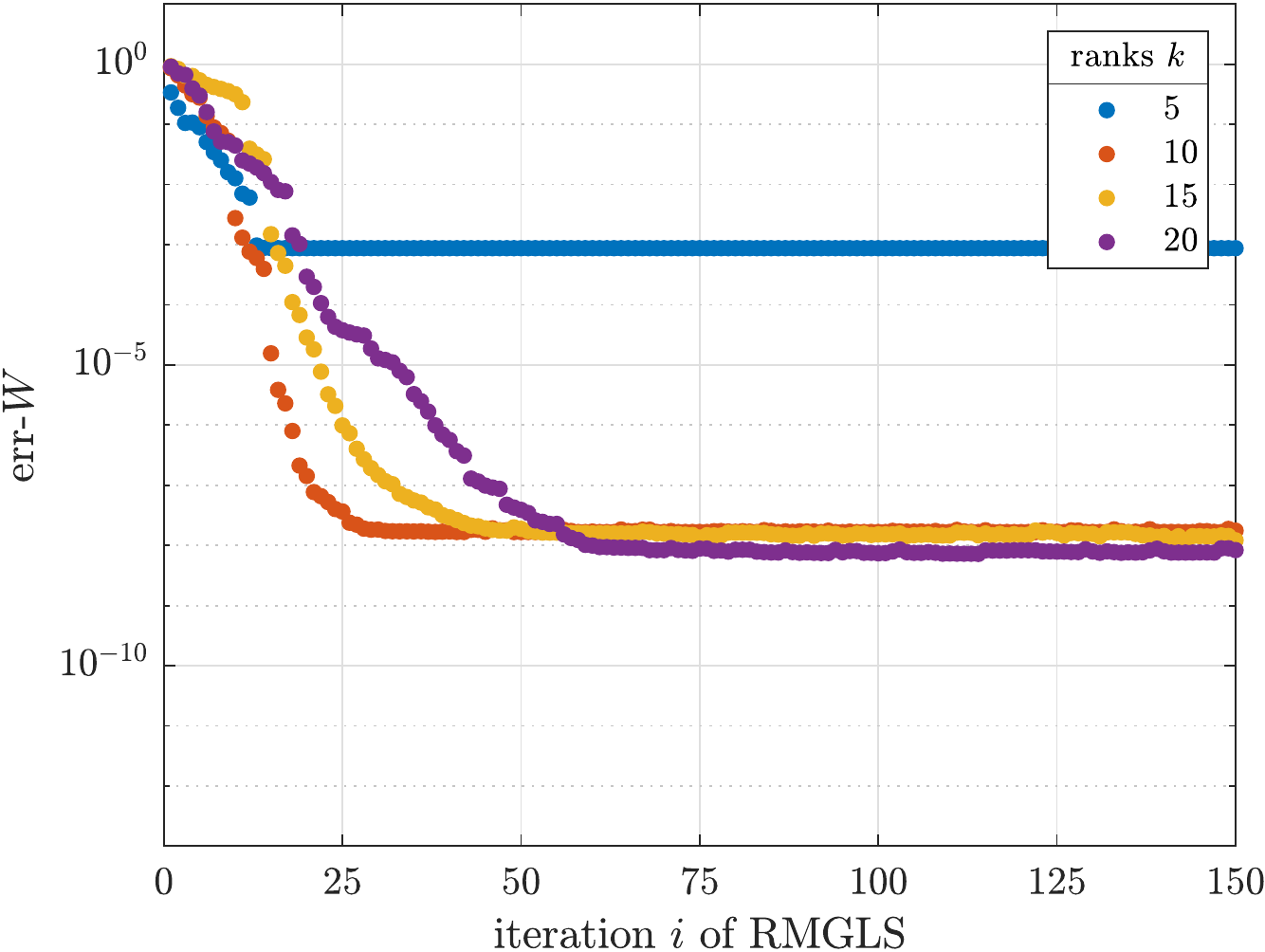}
        {\scriptsize (b) \emph{weak Wolfe}}
    \end{minipage}
    \caption{Convergence of $ \mathrm{err}\textrm{-}W $ for level $ \ell_{\mathrm{f}} = 8 $, for the problem of~\cref{sec:var_pb_1}.}\label{fig:LYAP_5_8_Conv_Uh_K5_K20_HZ_SS5_mgiter_150}
\end{figure}

\Cref{fig:LYAP_5_8_Conv_Uh_K5_K20_HZ_SS5_mgiter_150} shows the convergence of \textrm{err-$W$} for increasing ranks $k$. We compare a line search based on the new Hager--Zhang conditions to the weak Wolfe conditions. The plateaus in both panels are due to the fact that the approximate solution is computed in low-rank format and it is compared to the full-rank reference solution $ W^{(\ast)}_{h} $. The latter has good low-rank approximations, which is confirmed by the later onset of the stagnation phase when increasing the rank in RMGLS. Panel (b) of \cref{fig:LYAP_5_8_Conv_Uh_K5_K20_HZ_SS5_mgiter_150} shows that a line-search procedure with weak Wolfe conditions does not allow us to reduce \textrm{err-$W$} below $ \sqrt{\epsmach} \approx 10^{-8} $ in double precision arithmetic. This clearly makes the case that the Hager--Zhang line search is useful if we want to obtain more accurate low-rank approximations, as it is visible in panel (a).

To assess the accuracy of the solutions obtained for the Lyapunov equation, we also use the standard residual
\[
   r(W_h) \coloneqq \| AW_h + W_hA - \Gamma_h \|_{\F}.
\]
We also consider the following relative residual based on the backward error \protect{\cite[eq.~(3.6)]{Simoncini:2007}}
\[
   r_{\mathrm{BW}}(W_h) \coloneqq \frac{\| AW_h + W_hA - \Gamma_h \|_{\F}}{2\| A \|_{2} \| W_h \|_{\F} + \| \Gamma_h \|_{\F}}.
\]

\begin{figure}[htbp]
    \centering
    \begin{minipage}{0.48\textwidth}
        \centering
        \includegraphics[width=\textwidth]{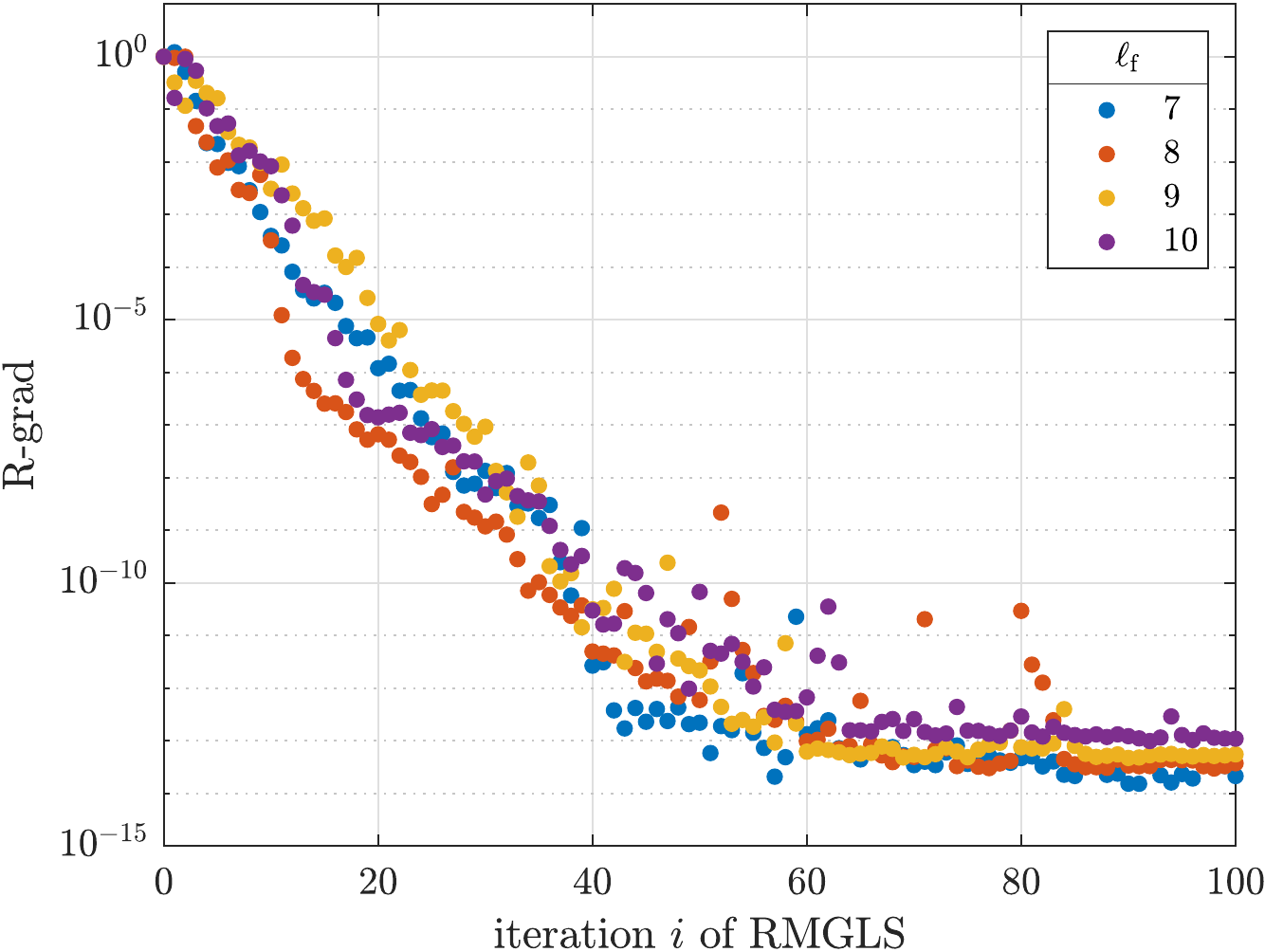}
          {\scriptsize (a) \emph{rank $k=5$}}        
    \end{minipage}\hfill
    \begin{minipage}{0.48\textwidth}
        \centering
        \includegraphics[width=\textwidth]{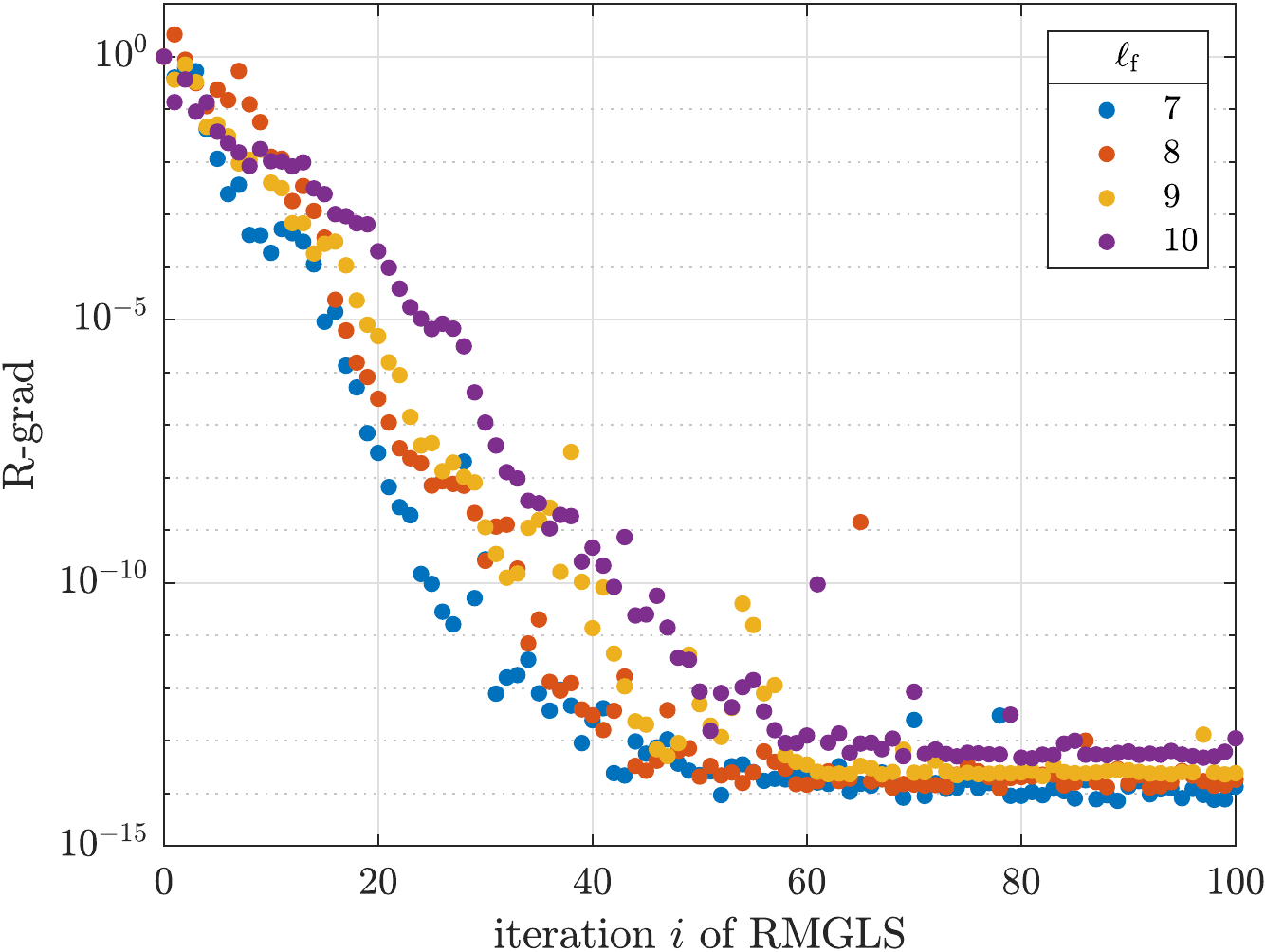}
          {\scriptsize (b) \emph{rank $k=10$}}        
    \end{minipage}
    \caption{Convergence of $ \mathrm{R\textrm{-}grad} $ for several finest levels $\ell_{\mathrm{f}}$, for the problem of \cref{sec:var_pb_1}.}\label{fig:LYAP_Conv_Grad_Vcycle_7_10_HZ_SS6_mgiter_100}        
\end{figure}

\Cref{fig:LYAP_Conv_Grad_Vcycle_7_10_HZ_SS6_mgiter_100} compares the convergence behavior of \textrm{R-grad} for different fine-scale manifolds with $ \ell_{\mathrm{f}} = 7, 8, 9, 10 $. The corresponding sizes of the discretizations are $16\,384$ (\textcolor[RGB]{23,126,194}{$ \bullet $}), $65\,536$ (\textcolor[RGB]{219,98,46}{$ \bullet $}), $262\,144$ (\textcolor[RGB]{237,183,52}{$ \bullet $}) and $1\,048\,576$ (\textcolor[RGB]{137,66,152}{$ \bullet $}). Panel (a) corresponds to rank $ k = 5 $, while panel (b) refers to $ k = 10 $.
One can observe that the convergence behavior is not very dependent on the mesh size, thereby confirming that RMGLS has an almost mesh-independent convergence typical of multigrid methods. 
In \cref{tab:var_pb_1_conv_grad_K5}, the final R-grad, err-$W$, backward error $ r_{\mathrm{BW}}(W_{h}) $ and residual $ r(W_{h}) $ of the Lyapunov equation are displayed.

\begin{table}[htbp]
   \caption{Final gradient norm and residuals for the problems of \cref{fig:LYAP_Conv_Grad_Vcycle_7_10_HZ_SS6_mgiter_100}. The error of the best rank-5 approximation is $ \approx 8.73 \times 10^{-4} $.  The error of the best rank-10 approximation is $ \approx 1.41 \times 10^{-8} $.} \label{tab:var_pb_1_conv_grad_K5}
   \vspace{-0.4cm}
\begin{center}
   \resizebox{\textwidth}{!}{%
      \begin{tabular}{c|cccccc}  
         \toprule
           &  $ \ell_{\mathrm{f}} $  &   size     &   R-grad(100)   &   $ \frac{r_{\mathrm{BW}}(W_{h}^{(100)})}{r_{\mathrm{BW}}(W_{h}^{(0)})} $   &  $ r(W_{h}^{(100)}) $   &   err-$W$(100)  \\
         \midrule
         \multirow{4}{*}{\STAB{\rotatebox[origin=c]{90}{rank 5}}}
             &  7 (\textcolor[RGB]{23,126,194}{$ \bullet $}) &       16\,384   &   $ 2.15 \times 10^{-14} $   &    $ 4.91 \times 10^{-5} $   &   $ 1.27 \times 10^{-4} $   &  $ 8.73 \times 10^{-4} $   \\
             &  8 (\textcolor[RGB]{219,98,46}{$ \bullet $})  &       65\,536   &   $ 3.76 \times 10^{-14} $   &    $ 1.65 \times 10^{-5} $    &   $ 6.34 \times 10^{-5} $   &  $ 8.74 \times 10^{-4} $   \\
             &  9 (\textcolor[RGB]{237,183,52}{$ \bullet $}) &      262\,144   &   $ 5.55 \times 10^{-14} $   &    $ 5.57 \times 10^{-6} $    &   $ 3.17 \times 10^{-5} $   &  $ 8.75 \times 10^{-4} $   \\
             & 10 (\textcolor[RGB]{137,66,152}{$ \bullet $}) &   1\,048\,576   &   $ 1.10 \times 10^{-13} $   &    $ 1.97 \times 10^{-6} $    &   $ 1.59 \times 10^{-5} $   &  $ 8.75 \times 10^{-4} $ \\
         \midrule
         \multirow{4}{*}{\STAB{\rotatebox[origin=c]{90}{rank 10}}}
             &  7 (\textcolor[RGB]{23,126,194}{$ \bullet $}) &       16\,384   &   $ 1.35 \times 10^{-14} $   &    $ 4.47 \times 10^{-9} $   &   $ 1.63 \times 10^{-8} $   &   $ 1.52 \times 10^{-8} $   \\
             &  8 (\textcolor[RGB]{219,98,46}{$ \bullet $})  &       65\,536   &   $ 1.83 \times 10^{-14} $   &    $ 1.54 \times 10^{-9} $    &   $ 8.46 \times 10^{-9} $   &   $ 1.54 \times 10^{-8} $   \\
             &  9 (\textcolor[RGB]{237,183,52}{$ \bullet $}) &      262\,144   &   $ 2.43 \times 10^{-14} $   &    $ 5.25 \times 10^{-10} $    &   $ 4.27 \times 10^{-9} $   &   $ 1.55 \times 10^{-8} $   \\
             & 10 (\textcolor[RGB]{137,66,152}{$ \bullet $}) &   1\,048\,576   &   $ 1.12 \times 10^{-13} $   &    $ 1.82 \times 10^{-10} $    &   $ 2.14 \times 10^{-9} $   &   $ 1.55 \times 10^{-8} $   \\
         \bottomrule
      \end{tabular}}
   \end{center}
\end{table}

The numerical experiments presented in this section show that RMGLS, our Riemannian multilevel optimization algorithm with Hager--Zhang line search, converges as we would expect from an effective multigrid method. Satisfying the approximate Wolfe conditions in the Hager--Zhang line search seems to be sufficient for the method to converge to local minima that are accurate when measured in the relative error and residual norms.

\subsubsection{Rank adaptivity}
\label{sec:rank_adaptivity}
\begin{figure}[htbp]
  \centering
  \includegraphics[width=0.65\textwidth]{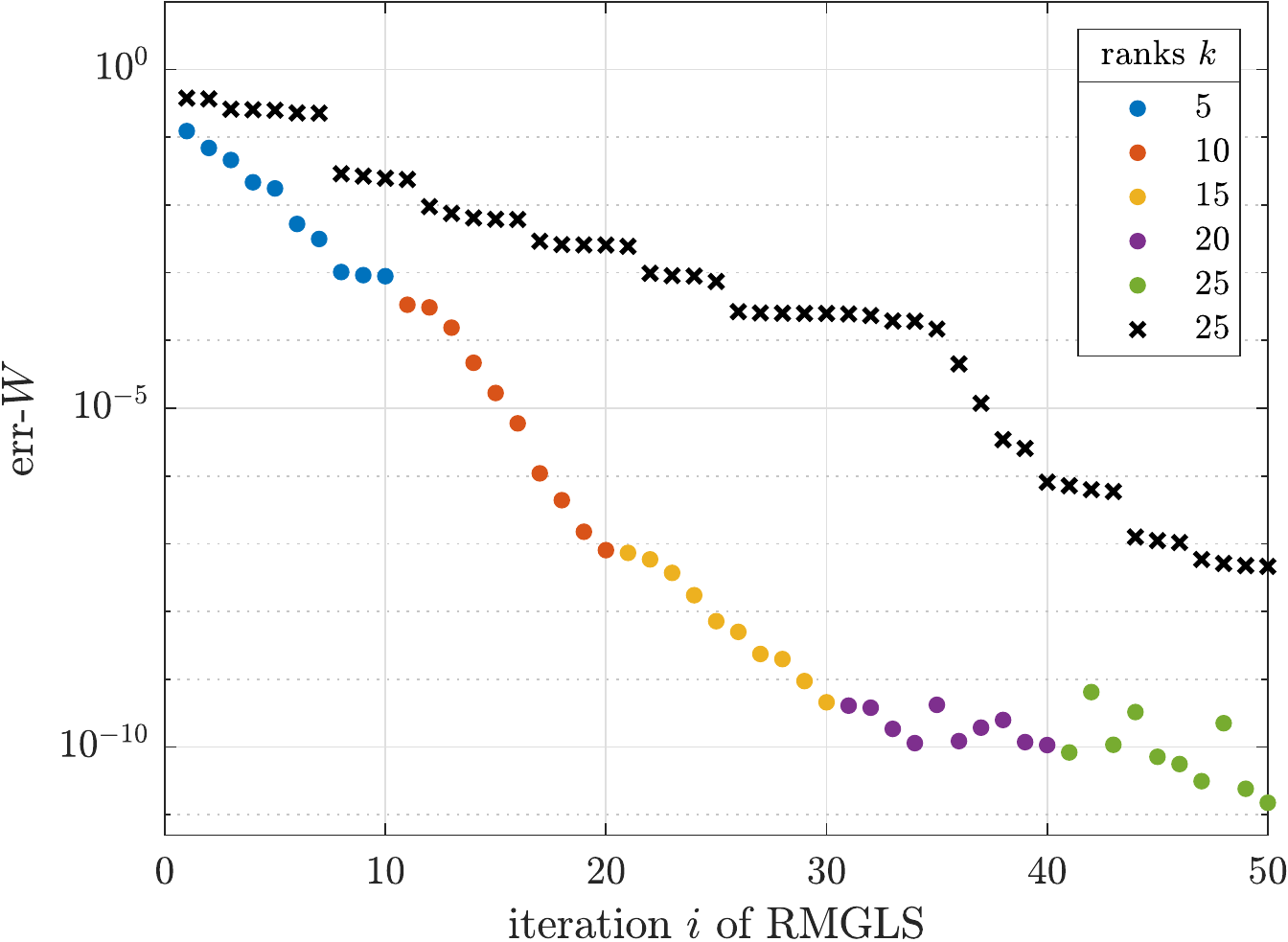}
  \caption{Rank-adaptivity for RMGLS applied to the problem of \cref{sec:var_pb_1}, with \cref{eq:rank5_rhs} as right-hand side. Starting from rank 5, the rank is increased by 5 every 10 iterations, until $ k = 25 $. The black crosses illustrate the behavior of non-adaptive RMGLS with rank $ k = 25 $.}
  \label{fig:LYAP_Conv_Uh_warm_start_5_25_HZ_SS5}
\end{figure}

In the framework of Riemannian optimization, rank-adaptivity can be introduced by successive runs of increasing rank, using the previous solution as a warm start for the next rank. For recent discussions about this approach see \protect{\cite{Uschmajew:2015,Kressner:2016}}. An example is given for the problem described in this section, with \cref{eq:rank5_rhs} as right-hand side, again with finest level $ \ell_{\mathrm{f}} = 8 $ and coarsest level $ \ell_{\mathrm{c}} = 5 $, using 5 smoothing steps.
Starting from rank $ k^{(0)} = 5 $, we run RMGLS for 10 iterations, and use the approximate solution to warm start the algorithm with ranks $ k^{(i)} = k^{(i-1)} + 5 $, $ i = 1, \ldots, 4 $. \Cref{fig:LYAP_Conv_Uh_warm_start_5_25_HZ_SS5} compares the convergence behavior of this adaptive strategy with the non-adaptive RMGLS, for a target rank $ k = 25 $. It is apparent that the adaptive RMGLS is more efficient than its non-adaptive counterpart. For example, at the 30th iteration, $ r(W_{h}^{(30)}) \approx 2.50 \times 10^{-4} $ for the non-adaptive RMGLS, whereas it is already $ r(W_{h}^{(30)}) \approx 4.57 \times 10^{-10} $ in the adaptive version.

\subsection{A nonlinear problem}
\label{sec:var_pb_2}
Next, we consider the variational problem from \protect{\cite[Example 5.1]{Wen:2009}} involving an exponential as nonlinear term:
\begin{equation}\label{eq:var_pb_2_functional}
   \begin{cases}
      \displaystyle\min_{w} \cF(w) = \int_{\Omega} \tfrac{1}{2} \| \nabla w \|^{2} + \lambda ( w - 1 ) \, e^{w} - \gamma \, w \dx\dy \\
      \quad \text{such that} \quad w=0 \ \text{on} \ \partial\Omega,
   \end{cases}
\end{equation}
where $\lambda = 10$, $\Omega = [0,1]^{2}$, and
\begin{equation*}
   \gamma(x,y) = \big( (9 \pi^2 + \lambda e^{(x^2-x^3)\sin(3\pi y)}) (x^2-x^3) + 6x - 2 \big) \sin(3\pi y).
\end{equation*}
The variational problem \cref{eq:var_pb_2_functional} corresponds to the nonlinear PDE~\protect{\cite[eq.~(5.4)]{Henson:2003}}
\begin{equation*}
\begin{cases}
   -\Delta w + \lambda w e^{w} - \gamma  = 0 \quad \text{in} \ \Omega, \\
   w = 0 \quad \text{on} \ \partial \Omega.
\end{cases}
\end{equation*}
The exact solution $ \wex = (x^2-x^3)\sin(3\pi y) $ has rank 1, making it less interesting as test case for our low-rank method. In addition, a discretization of the exponential term $ e^{w} $ does not admit a good low-rank approximation for $w$ close to the exact solution.

The following modification, 
\begin{equation}\label{eq:prob_var_2}
\begin{cases}
   -\Delta w + \lambda w (w+1) - \gamma  = 0 \quad \text{in} \ \Omega, \\
   w = 0 \quad \text{on} \ \partial \Omega,
\end{cases}
\end{equation}
is better suited as test case: as we will show below, the nonlinearity $w(w+1)$ can be computed efficiently when $w$ is low rank and the exact solution is full rank but has good low-rank approximations.

To obtain the variational problem corresponding to~\eqref{eq:prob_var_2}, $ -\Delta w $ gives rise to the term $ \frac{1}{2} \| \nabla w \|^{2} $ in the integrand of the objective functional, as seen in~\cref{sec:var_pb_1}. The term in $\gamma$ also remains the same. For the nonlinear term in the middle, we calculate the integral of $\lambda w (w+1) $ with respect to $w$, which gives $ \lambda w^{2} \big( \tfrac{1}{3}w + \tfrac{1}{2} \big) $.
Finally, we can formulate the variational problem as
\begin{equation}\label{eq:modified_NPDE_variational}
   \begin{cases}
      \displaystyle\min_{w} \cF(w) = \displaystyle\int_{\Omega} \tfrac{1}{2} \| \nabla w \|^{2} + \lambda w^{2} \big( \tfrac{1}{3}w + \tfrac{1}{2} \big) - \gamma \, w \dx\dy \\
      \quad \text{such that} \quad w=0 \ \text{on} \ \partial\Omega.
   \end{cases}
\end{equation}
For $ \gamma $, we choose the same right-hand side adopted in \cref{eq:rank5_rhs}.

\subsubsection{Discretization of the objective function}
Discretizing \cref{eq:modified_NPDE_variational} similarly as in \cref{sec:var_pb_1_discr_functional}, we obtain
\[
   \cF_{h} = h^{2} \sum_{i,j=0}^{2^{\ell} - 1} \left( \tfrac{1}{2} ( \partial w_{x_{ij}}^{2} + \partial w_{y_{ij}}^{2} ) + \lambda w_{ij}^{2}\big( \tfrac{1}{3}w_{ij} + \tfrac{1}{2} \big) - \gamma_{ij} w_{ij} \right).
\]
The first term and the third term have the same matrix form as the one seen in \cref{sec:var_pb_1_discr_functional}.
For the second term, we have
\[
   \sum_{i,j} \frac{\lambda}{2} \, w_{ij}^{2} = \frac{\lambda}{2} \trace(W_{h}\tr W_{h}) = \frac{\lambda}{2} \| \Sigma \|^{2}_{\F},
\]
and
\begin{equation}\label{eq:modified_third_term}
   \sum_{i,j} \frac{\lambda}{3} \, w_{ij}^{3} = \frac{\lambda}{3} \trace( W_{h}\tr ( W_{h} \odot W_{h} ) ).
\end{equation}
For the term $ W_{h} \odot W_{h} $, we perform the element-wise multiplication in factorized form as explained in \protect{\cite[§7]{Kressner:2014}} and store the result in the format $ U_{\odot}\Sigma_{\odot}V_{\odot}\tr $:
\[
   W_{h} \odot W_{h} = ( U \kt U ) ( \Sigma \otimes \Sigma ) (V \kt V)\tr = U_{\odot}\Sigma_{\odot}V_{\odot}\tr,
\]
where $ \kt $ denotes a transposed variant of the Khatri-Rao product (see definition in~\protect{\cite[§7]{Kressner:2014}}). Observe that $\rank(W_{h} \odot W_{h}) \leq k^2$.
As a consequence, the term \cref{eq:modified_third_term} becomes
\[
   \frac{\lambda}{3} \trace( W_{h}\tr ( W_{h} \odot W_{h} ) ) = \frac{\lambda}{3} \trace\!\big( V(U\Sigma)\tr U_{\odot}\Sigma_{\odot}V_{\odot}\tr \big) = \frac{\lambda}{3} \trace\!\big( \Sigma (U\tr U_{\odot}) \Sigma_{\odot} (V_{\odot}\tr V) \big).
\]
Finally, the discretized functional in matrix form is
\begin{align*}
    \cF_{h} = & \ \tfrac{h^{2}}{2} \Big( \| (LU) \Sigma \|^{2}_{\F} + \| (LV) \Sigma \|^{2}_{\F} + \lambda \| \Sigma \|^{2}_{\F} \\
    &  +  \ \tfrac{2}{3}\lambda \trace\!\big( \Sigma (U\tr U_{\odot}) \Sigma_{\odot}  (V_{\odot}\tr V) \big) - 2\trace\!\big( \Sigma_{\gamma} (U_{\gamma}\tr U)\Sigma (V\tr V_{\gamma}) \big) \Big),
\end{align*}
which can be evaluated in $ O\big(  nk (k_{\gamma} + k^2 ) + k ( k_{\gamma}^2 + k^3 )   \big) $ flops. 

\subsubsection{Discretization of the gradient}
The gradient of $ \cF $ is the functional derivative
\begin{equation*}
   \frac{\delta \cF}{\delta w} = -\Delta w + \lambda w (w+1) - \gamma.
\end{equation*}
The discretized Euclidean gradient in matrix form is given by
\[
    G_{h} = h^{2} \left( A W_{h} + W_{h} A + \lambda W_{h} \odot W_{h} + \lambda W_{h} - \Gamma_{h} \right),
\]
with $ A $ as in~\eqref{eq:var_pb_1_Laplacian_A}. 
Substituting the formats $ W_{h} = U\Sigma V\tr $, $ W_{h} \odot W_{h} = U_{\odot}\Sigma_{\odot}V_{\odot}\tr$, and $ \Gamma_{h} = U_{\gamma} \Sigma_{\gamma} V_{\gamma}\tr $, we get the factorized form
\begin{equation*}
  G_{h} = h^{2} \begin{bmatrix}
      (A+\lambda I)U  &  U  &  U_{\odot}  &  U_{\gamma}
   \end{bmatrix}
   \blkdiag\!\left(\Sigma, \Sigma, \lambda \Sigma_{\odot}, -\Sigma_{\gamma} \right)
   \begin{bmatrix}
      V  &  AV  &  V_{\odot}  &  V_{\gamma} 
   \end{bmatrix}\tr.
\end{equation*}
The gradient $G_h$ can be represented in only $O(nk)$ flops for computing $(A+\lambda I)U$ and $AV$.

\subsubsection{Numerical results}
We repeat the same set of experiments as for the previous problem to verify the convergence of the error and residual functions defined in~\cref{sec:num_res_var_pb_1}. The coarse level was again taken as $ \ell_{\mathrm{c}} = 5 $.

\begin{figure}[htbp]
    \centering
    \begin{minipage}{0.48\textwidth}
        \centering
        \includegraphics[width=\textwidth]{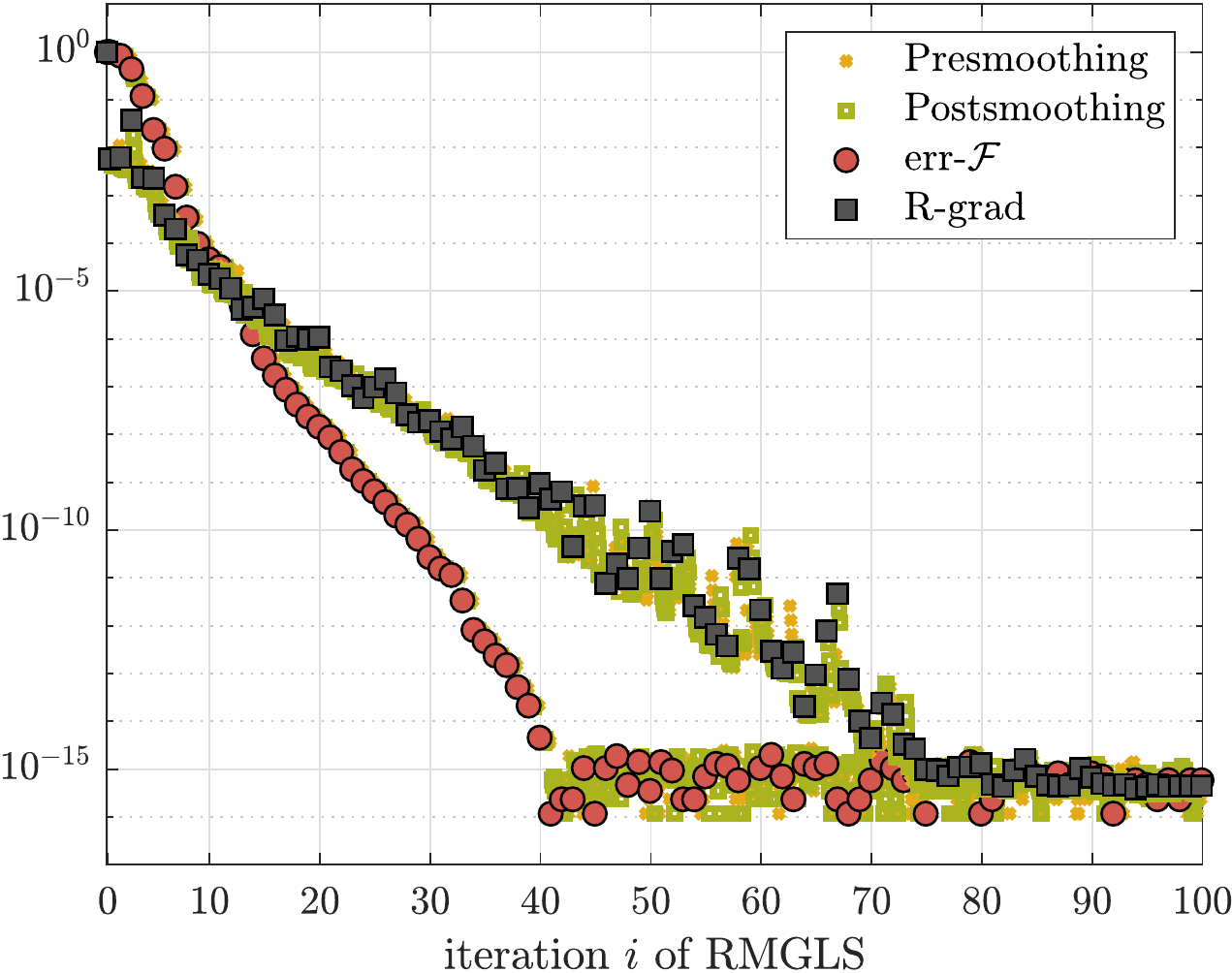}
        \caption{Convergence of $ \mathrm{err}\textrm{-}\cF $ and $ \mathrm{R\textrm{-}grad} $ for level $ \ell_{\mathrm{f}} = 8 $ and rank $ k = 5 $, for the problem of~\cref{sec:var_pb_2}.}\label{fig:NPDE_CostGrad_5_8_K5_HZ_SS5_mgiter_100}
    \end{minipage}\hfill
    \begin{minipage}{0.48\textwidth}
        \centering
        \includegraphics[width=\textwidth]{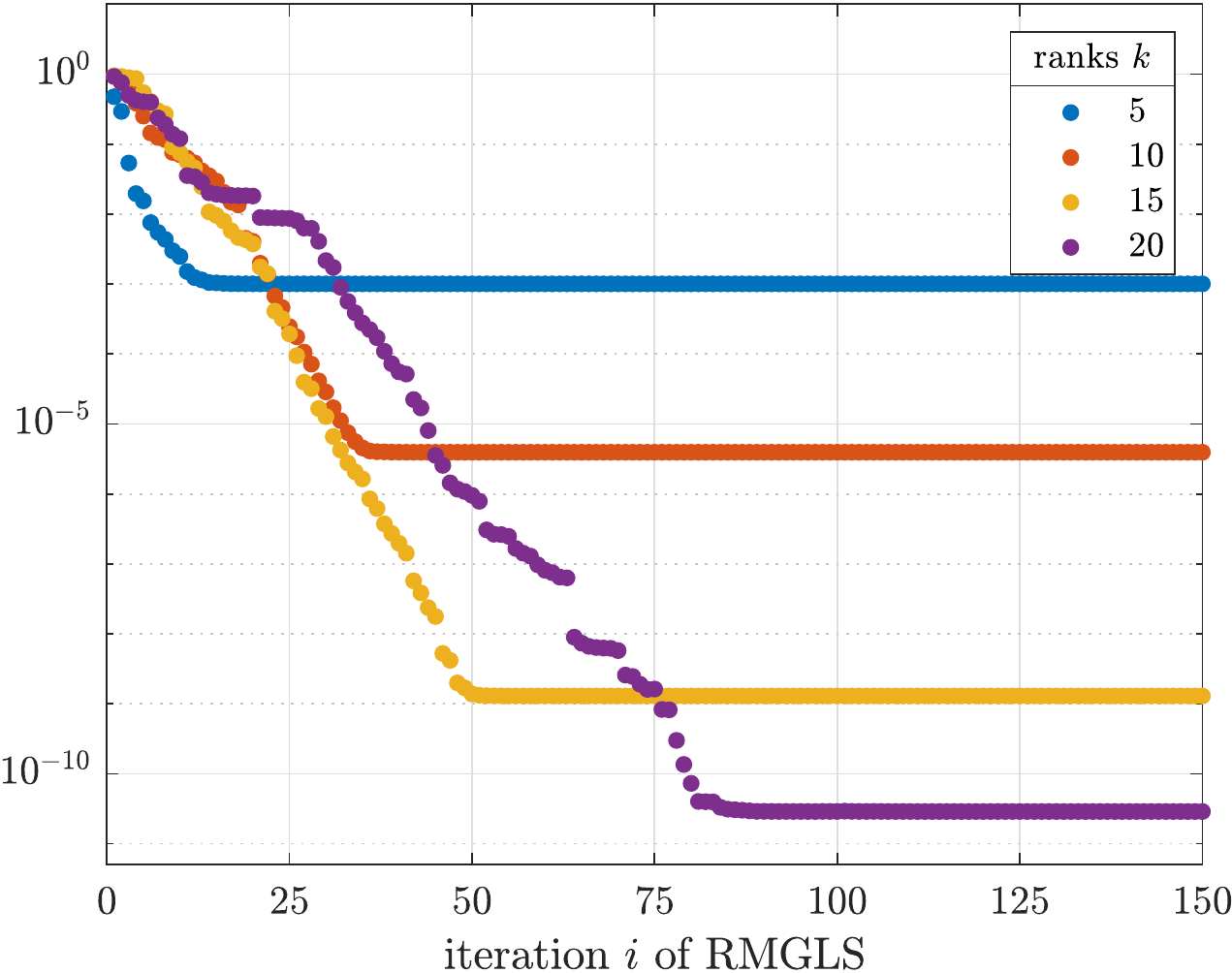}
        \caption{Convergence of $ \mathrm{err}\textrm{-}W $, with Hager--Zhang line search, for $ \ell_{\mathrm{f}} = 8 $ and the rank values $ k = 5,10,15,20 $.}\label{fig:NPDE_5_8_Conv_Uh_K5_K20_HZ_SS5_mgiter_150}
    \end{minipage}
\end{figure}

\begin{figure}[htbp]
  \centering
  \includegraphics[width=0.65\textwidth]{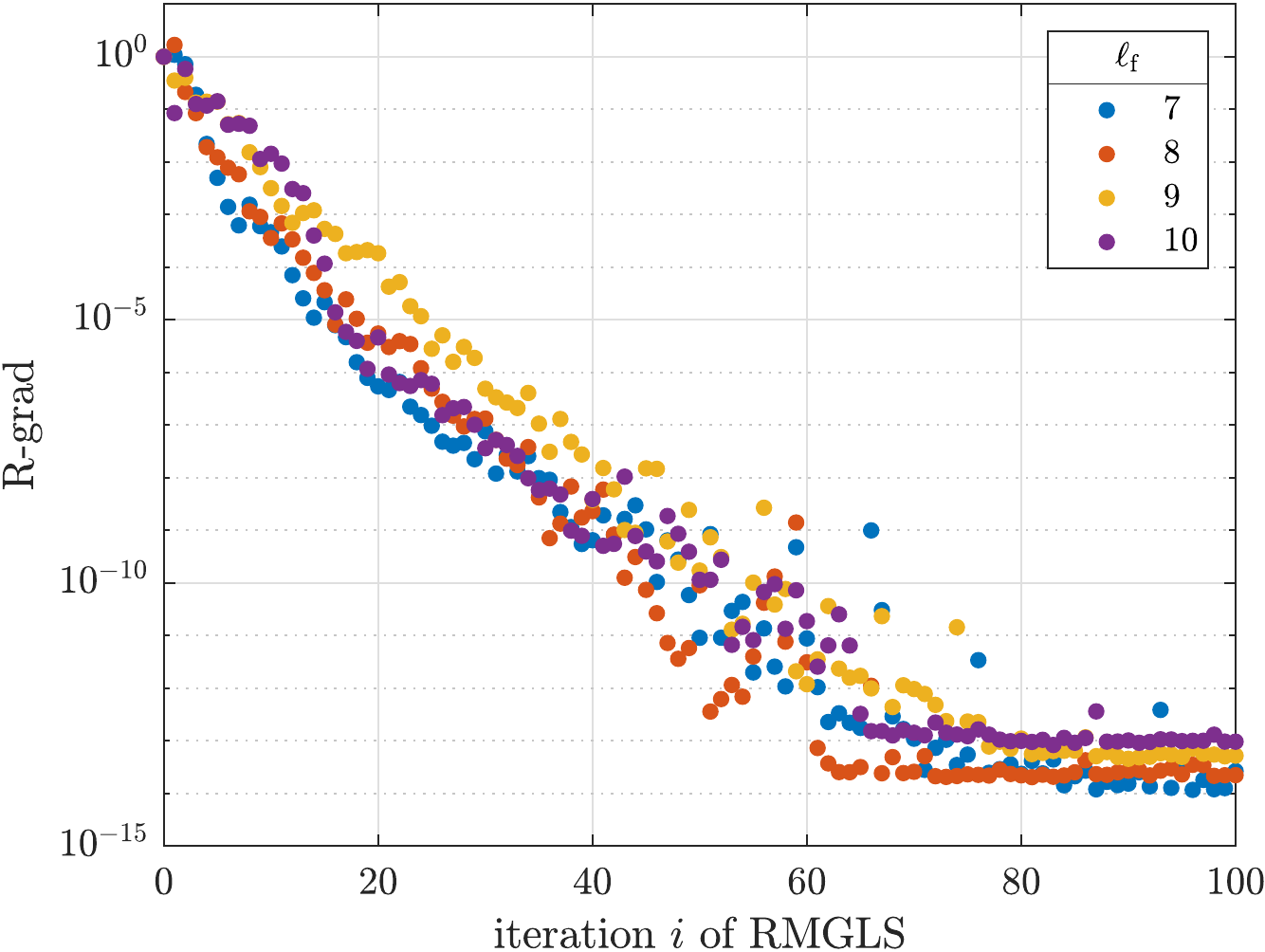}
  \caption{Convergence of $ \mathrm{R\textrm{-}grad} $ for several finest levels $\ell_{\mathrm{f}}$ and rank $ k = 5 $, for the problem of \cref{sec:var_pb_2}.}
  \label{fig:var_pb_2_conv_grad}
\end{figure}

\begin{figure}[htbp]
  \centering
  \includegraphics[width=0.65\textwidth]{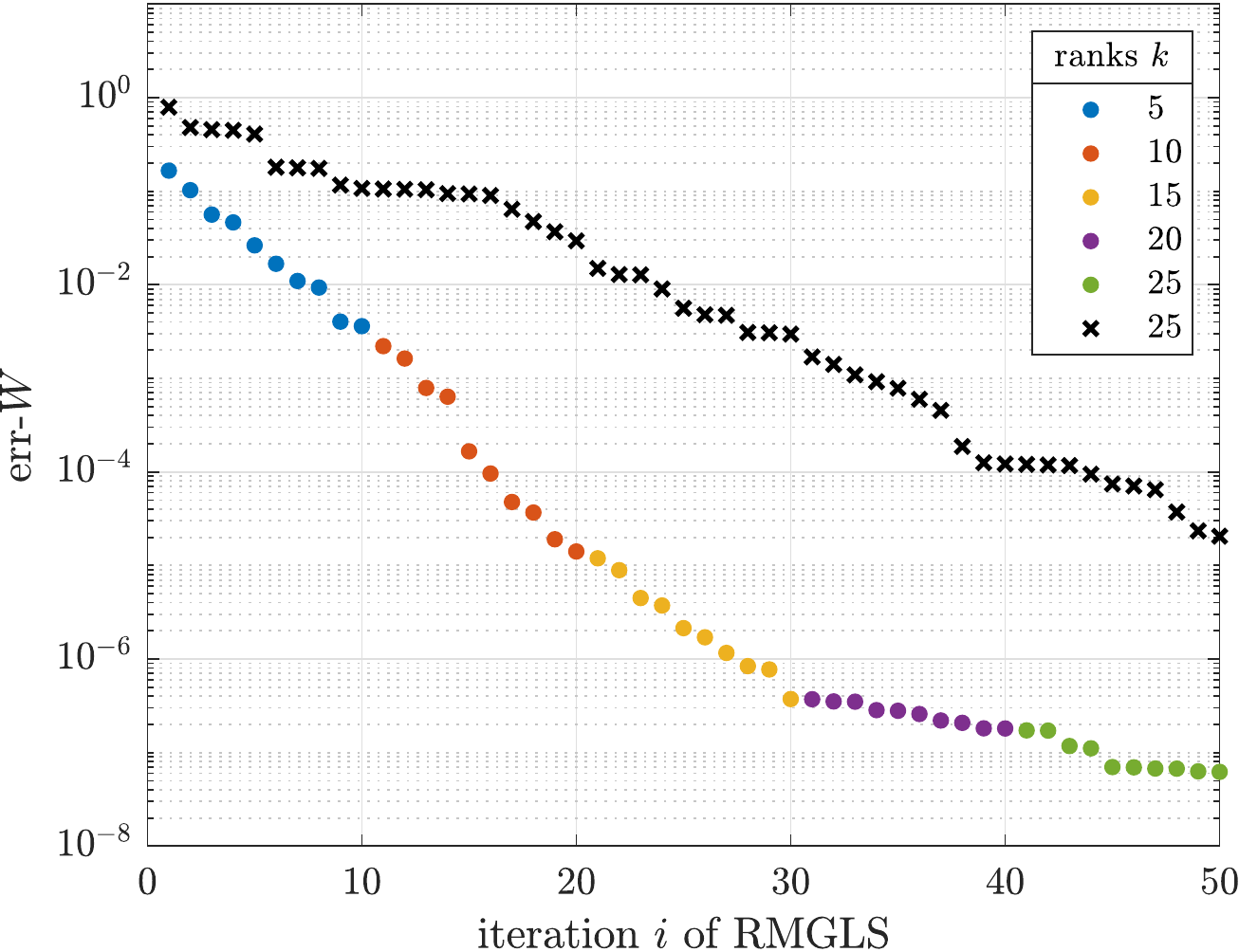}
  \caption{Rank-adaptivity for RMGLS applied to the problem of \cref{sec:var_pb_2}, with \cref{eq:rank5_rhs} as right-hand side. Starting from rank 5, the rank is increased by 5 every 10 iterations, until $ k = 25 $. The black crosses illustrate the behavior of non-adaptive RMGLS with rank $ k = 25 $.}
  \label{fig:NPDE_Conv_Uh_warm_start_5_25_HZ_SS5}
\end{figure}

Comparing \cref{fig:NPDE_CostGrad_5_8_K5_HZ_SS5_mgiter_100}, \ref{fig:NPDE_5_8_Conv_Uh_K5_K20_HZ_SS5_mgiter_150}, \ref{fig:var_pb_2_conv_grad}, and \ref{fig:NPDE_Conv_Uh_warm_start_5_25_HZ_SS5} for this nonlinear problem to the ones of the linear problem, we observe that the earlier conclusions remain virtually the same.

\section{Comparison with other methods}\label{sec:comparison}
We compare Euclidean trust-regions (ETR), Euclidean multilevel optimization (EML), EML with low rank via truncated SVD, Riemannian trust-regions (RTR) with fixed rank, and our RMGLS with fixed rank. ETR and EML do not use any low-rank approximation, whereas the other methods do. 

All methods were implemented in \textsc{Matlab}. ETR and RTR were executed using solvers from the Manopt package~\cite{Boumal:2014} with the Riemannian embedded submanifold geometry from~\cite{Vandereycken:2013} for RTR. EML was implemented by ourselves based on the same multigrid components as RMGLS, as already explained in~\cref{sec:variational_problems}. EML with truncated SVD applies truncation via the SVD with a fixed rank after every computational step in EML.

\Cref{tab:comparison_LYAP} summarizes the results for the linear problem described in \cref{sec:var_pb_1}. It is apparent that the Euclidean algorithms soon become very inefficient as the problem size grows. Hence we omit the results for bigger problem sizes. For the smaller problems, the residual of the final approximation was always very small since there was no rank truncation.

In the low-rank version of EML with truncated SVD, the algorithm was stopped \emph{before} stagnation in the residual norm started to occur, as determined by manual inspection. This was done so that the algorithm certainly did not run longer than needed.\footnote{Although in practice such a stopping condition can not be implemented.} All the other low-rank algorithms were stopped when the norm of the Riemannian gradient was below the threshold value of $ 10^{-12} $.

Observe that the accuracy achieved by EML with truncated SVD is not as good compared to RTR and RMGLS. This was \emph{not} due to our stopping condition but probably because of the fixed-rank truncations throughout the multigrid cycle in EML. It is possible that a more careful choice of ranks can improve on the accuracy, but we did not investigate this issue since RTR and RMGLS are also using fixed-rank truncations.

The Riemannian algorithms on the manifold of fixed-rank matrices show a more efficient behavior. For problems having a relatively small size ($ \ell_{\mathrm{f}} = 10, \ 11 $), RTR is more efficient than RMGLS, while for bigger problems, RMGLS is much more efficient that RTR. The fastest computational time for a given level is highlighted in bold text. In particular, for $ \ell_{\mathrm{f}} = 14 $, our RMGLS is almost 6 times more efficient than the RTR. This demonstrates that for very big problem sizes the Riemannian multilevel strategy is the most advantageous.

An important observation is that the Riemannian algorithms can be terminated based on the Riemannian gradient, since it provably can be made very small, as it is clear from the tables and also from the figures in the previous section. This property allows us to stop the algorithm when the local gradient is smaller than a certain threshold. On the contrary, the EML low-rank algorithm does not have this property and, since the (Riemannian) gradient might never become very small, the stopping criterion has to be based on stagnation detection.

Another observation concerns the ``multiplying factor'' across the levels for different methods. We are mostly interested in comparing the scaling factors for RTR and RMGLS when enlarging the level $\ell_{\mathrm{f}}$, since the other methods are visibly more expensive than these two. From \cref{tab:comparison_LYAP}, we obtain on average the scaling factors of 3.5 for RTR and 1.7 for RMGLS, respectively.

Finally, \cref{tab:comparison_EML_ranks_LYAP} shows that if we increase the rank, then it is possible to achieve better accuracy in the residuals $ r(W_{h}^{(\mathrm{end})}) $ for both EML with rank truncation and RMGLS. In addition, RMGLS is considerably faster than EML for the same rank and for the biggest problem.

\Cref{tab:comparison_NPDE} summarizes the results for the nonlinear problem described in \cref{sec:var_pb_2}. Similar considerations as above can be done for this problem. We point out that the higher computational times in the low-rank algorithms are due to the calculations of the Hadamard products in factored form. From \cref{tab:comparison_NPDE} we can obtain the following average multiplying factors across the levels: 4.3 for RTR, 2.0 for RMGLS. These are in good agreement with the ones computed for the linear problem.

\begin{table}[p]
   \caption{Comparisons of different methods for the problem described in \cref{sec:var_pb_1}. 
   The --- means the Riemannian gradient $\xi_{h}^{(\mathrm{end})}$ does not apply.}
   \label{tab:comparison_LYAP}
   \begin{center}
      \begin{tabular}{crccc}  
         \toprule
          $ \ell_{\mathrm{f}} $  &   size     &  time (s)  &  $ \| \xi_{h}^{(\mathrm{end})} \|_{F} $  &   $ r(W_{h}^{(\mathrm{end})}) $   \\
         \midrule \midrule
         \multicolumn{5}{c}{ETR (no rank truncation, no multilevel)}  \\
         \midrule
             9  &     262\,144  &     19 & ---  &  $ 9.2451 \times 10^{-15} $  \\
            10  &  1\,048\,576  &    164 & ---  &  $ 5.2284 \times 10^{-15} $  \\
            11  &  4\,194\,304  & 1\,787 & ---  &  $ 1.0223 \times 10^{-14} $  \\
         \midrule \midrule
         \multicolumn{5}{c}{EML (8 smoothing steps, $ \ell_{\mathrm{c}} = 7 $)}  \\
         \midrule
         \multicolumn{5}{c}{no rank truncation}  \\
         \midrule
             9  &     262\,144  &     16 & ---  &  $ 6.2645 \times 10^{-13} $  \\
            10  &  1\,048\,576  &     77 & ---  &  $ 3.4368 \times 10^{-13} $  \\
            11  &  4\,194\,304  &    459 & ---  &  $ 4.2710 \times 10^{-12} $  \\
         \midrule
         \multicolumn{5}{c}{truncation to rank 5}  \\
         \midrule
             9  &     262\,144  &     9  & --- &  $ 4.5166 \times 10^{-5} $  \\
            10  &  1\,048\,576  &    35  & --- &  $ 2.2084 \times 10^{-5} $  \\
            11  &  4\,194\,304  &    58  & --- &  $ 1.7780 \times 10^{-4} $  \\  
         \midrule \midrule
         \multicolumn{5}{c}{RTR -- rank 5 (no multilevel)}  \\
         \midrule
            10  &   1\,048\,576 &   \textbf{6} &  $ 1.5002 \times 10^{-14} $  &  $ 1.5873 \times 10^{-5} $  \\
            11  &   4\,194\,304 &  \textbf{20} &  $ 2.7687 \times 10^{-14} $  &  $ 7.9369 \times 10^{-6} $  \\
            12  &  16\,777\,216 &          66  &  $ 6.6810 \times 10^{-14} $  &  $ 3.9685 \times 10^{-6} $  \\
            13  &  67\,108\,864 &         237  &  $ 1.1654 \times 10^{-13} $  &  $ 1.9842 \times 10^{-6} $  \\
            14  & 268\,435\,456 &         929  &  $ 2.6852 \times 10^{-13} $  &  $ 9.9212 \times 10^{-7} $  \\
         \midrule \midrule
         \multicolumn{5}{c}{RMGLS -- rank 5 (8 smoothing steps, $ \ell_{\mathrm{c}} = 7 $)}  \\
         \midrule 
            10  &   1\,048\,576 &           18  &  $ 6.1634 \times 10^{-13} $  &  $ 1.5873 \times 10^{-5} $  \\
            11  &   4\,194\,304 &           26  &  $ 2.5091 \times 10^{-13} $  &  $ 7.9369 \times 10^{-6} $  \\
            12  &  16\,777\,216 &   \textbf{52} &  $ 6.5807 \times 10^{-13} $  &  $ 3.9685 \times 10^{-6} $  \\
            13  &  67\,108\,864 &   \textbf{94} &  $ 9.2574 \times 10^{-13} $  &  $ 1.9842 \times 10^{-6} $  \\
            14  & 268\,435\,456 &  \textbf{161} &  $ 6.1323 \times 10^{-13} $  &  $ 9.9212 \times 10^{-7} $  \\ 
         \bottomrule
      \end{tabular}
   \end{center}
\end{table}

\begin{table}[htbp]
   \caption{Comparisons of EML with rank truncation and RMGLS for different ranks applied to the problem described in \cref{sec:var_pb_1}. In both cases, 8 smoothing steps and coarsest level $ \ell_{\mathrm{c}} = 7 $ are used.} \label{tab:comparison_EML_ranks_LYAP}\vspace{-0.3cm}
      \begin{center}
   \resizebox{\textwidth}{!}{%
      \begin{tabular}{l|crccccc}  
         \toprule
            \multicolumn{1}{c}{}  &  &  & \multicolumn{2}{c}{EML}  &  \multicolumn{3}{c}{RMGLS} \\
         \cmidrule(lr){4-5} \cmidrule(lr){6-8}

            \multicolumn{1}{c}{}  &  $ \ell_{\mathrm{f}} $  &   size     &  time (s)  &  $ r(W_{h}^{(\mathrm{end})}) $ &  time (s)  &  $ \| \xi_{h}^{(\mathrm{end})} \|_{F} $   &  $ r(W_{h}^{(\mathrm{end})}) $  \\
         \midrule
         \multirow{3}{*}{\STAB{\rotatebox[origin=c]{90}{rank 10}}}
             &  9  &     262\,144  &    \textbf{13}  &  $ 9.1637 \times 10^{-9} $  &  18  &  $ 7.6740 \times 10^{-13} $  &  $ 4.2704 \times 10^{-9} $  \\
             & 10  &  1\,048\,576  &    55  &  $ 3.3303 \times 10^{-9} $  &  \textbf{31}  &  $ 5.0568 \times 10^{-13} $  &  $ 2.1431 \times 10^{-9} $  \\
             & 11  &  4\,194\,304  &   451  &  $ 4.9866 \times 10^{-5} $  &  \textbf{76}  &  $ 3.1722 \times 10^{-13} $  &  $ 1.0704 \times 10^{-9} $  \\
         \midrule
         \multirow{3}{*}{\STAB{\rotatebox[origin=c]{90}{rank 15}}}
             &  9  &     262\,144  &    43  &  $ 4.7685 \times 10^{-11} $  &  \textbf{42}  &  $ 6.1597 \times 10^{-13} $  &  $ 4.2953 \times 10^{-11} $  \\
             & 10  &  1\,048\,576  &   107  &  $ 2.4681 \times 10^{-11} $  & \textbf{103}  &  $ 6.1486 \times 10^{-13} $  &  $ 2.1541 \times 10^{-11} $  \\
             & 11  &  4\,194\,304  &   495  &  $ 4.5356 \times 10^{-11} $  &  \textbf{174}  &  $ 8.9919 \times 10^{-13} $  &  $ 1.0940 \times 10^{-11} $  \\
         \bottomrule
      \end{tabular}}
   \end{center}
\end{table}

\begin{table}[htbp]
   \caption{Comparisons of different methods for the problem described in \cref{sec:var_pb_2}. 
   The --- means the Riemannian gradient $\xi_{h}^{(\mathrm{end})}$ does not apply.}
   \label{tab:comparison_NPDE}
   \begin{center}
      \begin{tabular}{crccc}  
         \toprule
         $ \ell_{\mathrm{f}} $  &   size     &  time (s)  &  $ \| \xi_{h}^{(\mathrm{end})} \|_{F} $  &  $ r(W_{h}^{(\mathrm{end})}) $  \\
         \midrule \midrule
         \multicolumn{5}{c}{ETR: no rank truncation, no multilevel}  \\
         \midrule
             9  &     262\,144  &       23  & ---  &  $ 8.8890 \times 10^{-15} $  \\
            10  &  1\,048\,576  &      196  & ---  &  $ 6.0243 \times 10^{-15} $  \\
            11  &  4\,194\,304  &   1\,990  & ---  &  $ 1.1352 \times 10^{-14} $  \\
         \midrule \midrule
         \multicolumn{5}{c}{EML: $ \ell_{\mathrm{c}} = 7 $, 8 smoothing steps}  \\
         \midrule
         \multicolumn{5}{c}{no rank truncation}  \\
         \midrule
             9  &     262\,144  &    22  & ---  &  $ 2.6912 \times 10^{-13} $  \\
            10  &  1\,048\,576  &    75  & ---  &  $ 9.4184 \times 10^{-13} $  \\
            11  &  4\,194\,304  &   506  & ---  &  $ 4.2292 \times 10^{-13} $  \\
         \midrule
         \multicolumn{5}{c}{truncation to rank 5}  \\
         \midrule
             9  &     262\,144  &    11  & ---  &  $ 4.4994 \times 10^{-5} $  \\
            10  &  1\,048\,576  &    47  & ---  &  $ 5.0028 \times 10^{-5} $  \\
            11  &  4\,194\,304  &    72  & ---  &  $ 1.6216 \times 10^{-4} $  \\
         \midrule \midrule
         \multicolumn{5}{c}{RTR: rank 5, no multilevel}  \\
         \midrule
            10  &   1\,048\,576 &   \textbf{11} &  $ 1.3449 \times 10^{-14} $  &  $ 1.5614 \times 10^{-5} $  \\
            11  &   4\,194\,304 &   \textbf{31} &  $ 9.3240 \times 10^{-14} $  &  $ 7.8072 \times 10^{-6} $  \\
            12  &  16\,777\,216 &          151  &  $ 5.9424 \times 10^{-14} $  &  $ 3.9036 \times 10^{-6} $  \\
            13  &  67\,108\,864 &          554  &  $ 1.1696 \times 10^{-13} $  &  $ 1.9518 \times 10^{-6} $  \\
            14  & 268\,435\,456 &       3\,338  &  $ 2.1950 \times 10^{-13} $  &  $ 9.7591 \times 10^{-7} $  \\
         \midrule \midrule
         \multicolumn{5}{c}{RMGLS: rank 5, $ \ell_{\mathrm{c}} = 7 $, 8 smoothing steps}  \\
         \midrule 
            10  &   1\,048\,576 &           44  &  $ 5.3255 \times 10^{-13} $  &  $ 1.5614 \times 10^{-5} $  \\
            11  &   4\,194\,304 &           51  &  $ 4.0362 \times 10^{-13} $  &  $ 7.8072 \times 10^{-6} $  \\
            12  &  16\,777\,216 &  \textbf{120} &  $ 9.6698 \times 10^{-13} $  &  $ 3.9036 \times 10^{-6} $  \\
            13  &  67\,108\,864 &  \textbf{209} &  $ 3.8296 \times 10^{-13} $  &  $ 1.9518 \times 10^{-6} $  \\
            14  & 268\,435\,456 &  \textbf{549} &  $ 9.8448 \times 10^{-13} $  &  $ 9.7591 \times 10^{-7} $  \\
         \bottomrule
      \end{tabular}
   \end{center}
\end{table}

\begin{table}[htbp]
   \caption{Comparisons of EML with rank truncation and RMGLS for different ranks applied to the problem described in \cref{sec:var_pb_2}. In both cases, 8 smoothing steps and coarsest level $ \ell_{\mathrm{c}} = 7 $ are used.} \label{tab:comparison_EML_ranks_NPDE}\vspace{-0.3cm}
      \begin{center}
   \resizebox{\textwidth}{!}{%
      \begin{tabular}{l|crccccc}  
         \toprule
            \multicolumn{1}{c}{}  &  &  & \multicolumn{2}{c}{EML}  &  \multicolumn{3}{c}{RMGLS} \\
         \cmidrule(lr){4-5} \cmidrule(lr){6-8} 

            \multicolumn{1}{c}{}  &  $ \ell_{\mathrm{f}} $  &   size     &  time (s)  &  $ r(W_{h}^{(\mathrm{end})}) $ &  time (s)  &  $ \| \xi_{h}^{(\mathrm{end})} \|_{F} $   &  $ r(W_{h}^{(\mathrm{end})}) $  \\
         \midrule
         \multirow{3}{*}{\STAB{\rotatebox[origin=c]{90}{rank 10}}}
            &  9  &     262\,144  &    30  &  $ 4.7324 \times 10^{-7} $  &  \textbf{21}  &  $ 7.8437 \times 10^{-13} $  &  $ 3.7321 \times 10^{-7} $  \\
            & 10  &  1\,048\,576  &   123  &  $ 3.4975 \times 10^{-7} $  &  \textbf{61}  &  $ 4.0398 \times 10^{-13} $  &  $ 1.8660 \times 10^{-7} $  \\
            & 11  &  4\,194\,304  &   797  &  $ 1.2826 \times 10^{-5} $  &  \textbf{153}  &  $ 5.5800 \times 10^{-13} $  &  $ 9.3301 \times 10^{-8} $  \\
         \midrule
         \multirow{3}{*}{\STAB{\rotatebox[origin=c]{90}{rank 15}}}
            &  9  &     262\,144  &    107 &  $ 7.4928 \times 10^{-10} $  &  \textbf{92}   &  $ 2.0183 \times 10^{-13} $  &  $ 4.2886 \times 10^{-10} $  \\
            & 10  &  1\,048\,576  &    380 &  $ 9.6225 \times 10^{-10} $  &  \textbf{207}  &  $ 6.5306 \times 10^{-13} $  &  $ 2.6044 \times 10^{-10} $  \\
            & 11  &  4\,194\,304  & 3\,113 &  $ 4.3682 \times 10^{-10} $  &  \textbf{532}  &  $ 1.3610 \times 10^{-13} $  &  $ 8.3563 \times 10^{-11} $  \\
         \bottomrule
      \end{tabular}}
   \end{center}
\end{table}

\section{Conclusions} \label{sec:conclusions}

In this paper, we have shown how to combine multilevel optimization with optimization on low-rank manifolds. Compared to other approaches that are based on alternating minimization or rank-truncated iterative methods, our method does not need to accurately solve ill-conditioned local systems or use explicit preconditioning thanks to the intrinsic scalability of multigrid algorithms. This is demonstrated in our numerical experiments for two variational problems where our method succeeds in computing good low-rank approximations with an almost mesh-independent convergence behavior. In addition, we discussed an accurate line-search method to obtain highly accurate stationary points when only using first-order gradient information.

\bibliographystyle{siamplain}
\bibliography{references}

\begin{thebibliography}{10}

\bibitem{AMS:2008}
{\sc P.-A. Absil, R.~Mahony, and R.~Sepulchre}, {\em {Optimization Algorithms
  on Matrix Manifolds}}, Princeton University Press, Princeton, NJ, 2008.

\bibitem{Absil:2012}
{\sc P.-A. Absil and J.~Malick}, {\em Projection-like retractions on matrix
  manifolds}, SIAM Journal on Optimization, 22 (2012), pp.~135--158,
  \url{https://doi.org/10.1137/100802529}.

\bibitem{Absil:2015}
{\sc P.-A. Absil and I.~V. Oseledets}, {\em Low-rank retractions: a survey and
  new results}, Computational Optimization and Applications, 62 (2015),
  pp.~5--29, \url{https://doi.org/10.1007/s10589-014-9714-4}.

\bibitem{Boumal:2014}
{\sc N.~Boumal, B.~Mishra, P.-A. Absil, and R.~Sepulchre}, {\em {M}anopt, a
  {M}atlab toolbox for optimization on manifolds}, Journal of Machine Learning
  Research, 15 (2014), pp.~1455--1459, \url{http://www.manopt.org}.

\bibitem{Brenner:2007}
{\sc S.~Brenner and R.~Scott}, {\em The Mathematical Theory of Finite Element
  Methods}, Texts in Applied Mathematics, Springer New York, 2007.

\bibitem{Elman:2018}
{\sc H.~C. Elman and T.~Su}, {\em {A Low-Rank Multigrid Method for the
  Stochastic Steady-State Diffusion Problem}}, SIAM Journal on Matrix Analysis
  and Applications, 39 (2018), pp.~492--509,
  \url{https://doi.org/10.1137/17M1125170}.

\bibitem{Grasedyck:2007}
{\sc L.~Grasedyck and W.~Hackbusch}, {\em {A Multigrid Method to Solve Large
  Scale Sylvester Equations}}, SIAM Journal on Matrix Analysis and
  Applications, 29 (2007), pp.~870--894,
  \url{https://doi.org/10.1137/040618102}.

\bibitem{Gratton:2010}
{\sc S.~Gratton, M.~Mouffe, A.~Sartenaer, P.~L. Toint, and D.~Tomanos}, {\em
  Numerical experience with a recursive trust-region method for multilevel
  nonlinear bound-constrained optimization}, Optimization Methods and Software,
  25 (2010), pp.~359--386, \url{https://doi.org/10.1080/10556780903239295}.

\bibitem{Gratton:2008}
{\sc S.~Gratton, A.~Sartenaer, and P.~L. Toint}, {\em Recursive trust-region
  methods for multiscale nonlinear optimization}, SIAM Journal on Optimization,
  19 (2008), pp.~414--444, \url{https://doi.org/10.1137/050623012}.

\bibitem{Hackbusch:1985}
{\sc W.~Hackbusch}, {\em Multi-grid methods and applications}, Springer, 1985.

\bibitem{Hackbusch:2012}
{\sc W.~Hackbusch}, {\em Tensor Spaces and Numerical Tensor Calculus},
  Springer, 2012.

\bibitem{Hager:2005}
{\sc W.~W. Hager and H.~Zhang}, {\em A new conjugate gradient method with
  guaranteed descent and an efficient line search}, SIAM Journal on
  Optimization, 16 (2005), pp.~170--192,
  \url{https://doi.org/10.1137/030601880}.

\bibitem{Hager:2006}
{\sc W.~W. Hager and H.~Zhang}, {\em {Algorithm 851: CG\_DESCENT, a Conjugate
  Gradient Method with Guaranteed Descent}}, ACM Trans. Math. Softw., 32
  (2006), pp.~113--137, \url{https://doi.org/10.1145/1132973.1132979}.

\bibitem{Henson:2003}
{\sc V.~E. Henson}, {\em {Multigrid methods nonlinear problems: an overview}},
  in Computational Imaging, C.~A. Bouman and R.~L. Stevenson, eds., vol.~5016,
  International Society for Optics and Photonics, SPIE, 2003, pp.~36--48,
  \url{https://doi.org/10.1117/12.499473}.

\bibitem{Kressner:2016}
{\sc D.~Kressner, M.~Steinlechner, and B.~Vandereycken}, {\em Preconditioned
  low-rank riemannian optimization for linear systems with tensor product
  structure}, SIAM Journal on Scientific Computing, 38 (2016),
  pp.~A2018--A2044, \url{https://doi.org/10.1137/15M1032909}.

\bibitem{Kressner:2014}
{\sc D.~Kressner and C.~Tobler}, {\em {Algorithm 941: Htucker---A Matlab
  Toolbox for Tensors in Hierarchical Tucker Format}}, ACM Trans. Math. Softw.,
  40 (2014), pp.~22:1--22:22, \url{https://doi.org/10.1145/2538688}.

\bibitem{LeDret:2016}
{\sc H.~Le~Dret and B.~Lucquin}, {\em {Partial Differential Equations:
  Modeling, Analysis and Numerical Approximation}}, Birkhäuser, Basel, 2016.

\bibitem{Lewis:2005}
{\sc R.~M. Lewis and S.~G. Nash}, {\em Model problems for the multigrid
  optimization of systems governed by differential equations}, SIAM J. Sci.
  Comput., 26 (2005), pp.~1811--1837,
  \url{https://doi.org/10.1137/S1064827502407792}.

\bibitem{Luenberger:1972}
{\sc D.~G. Luenberger}, {\em The gradient projection method along geodesics},
  Management Science, 18 (1972), pp.~620--631,
  \url{http://www.jstor.org/stable/2629156}.

\bibitem{Mishra:2011b}
{\sc B.~Mishra, G.~Meyer, F.~Bach, and R.~Sepulchre}, {\em Low-rank
  optimization with trace norm penalty}, SIAM Journal on Optimization, 23
  (2013), pp.~2124--2149.

\bibitem{Mishra_V_2014}
{\sc B.~Mishra and B.~Vandereycken}, {\em A {R}iemannian approach to low-rank
  algebraic {R}iccati equations}, in 21st International Symposium on
  Mathematical Theory of Networks and Systems, 2014.

\bibitem{Nash:2000}
{\sc S.~G. Nash}, {\em A multigrid approach to discretized optimization
  problems}, Optimization Methods and Software, 14 (2000), pp.~99--116,
  \url{https://doi.org/10.1080/10556780008805795}.

\bibitem{Penzl:1997}
{\sc T.~Penzl}, {\em {A Multi-Grid Method for Generalized Lyapunov Equations}},
  {Technical Report SFB393/97--24}, Fakultät für Mathematik, TU Chemnitz,
  Chemnitz, Germany, 1997, \url{http://www.tu-chemnitz.de/sfb393/sfb97pr.html}.

\bibitem{RAKHUBA2019718}
{\sc M.~Rakhuba, A.~Novikov, and I.~Oseledets}, {\em Low-rank {R}iemannian
  eigensolver for high-dimensional {H}amiltonians}, Journal of Computational
  Physics, 396 (2019), pp.~718--737.

\bibitem{Rakhuba:2018ab}
{\sc M.~Rakhuba and I.~Oseledets}, {\em Jacobi--{Davidson} {Method} on
  {Low}-{Rank} {Matrix} {Manifolds}}, SIAM Journal on Scientific Computing, 40
  (2018), pp.~A1149--A1170.

\bibitem{Rosen:1995}
{\sc I.~G. Rosen and C.~Wang}, {\em {A Multilevel Technique for the Approximate
  Solution of Operator Lyapunov and Algebraic Riccati Equations}}, SIAM Journal
  on Numerical Analysis, 32 (1995), pp.~514--541,
  \url{http://www.jstor.org/stable/2158409}.

\bibitem{Rosen:1961}
{\sc J.~B. Rosen}, {\em The gradient projection method for nonlinear
  programming. {P}art {II}. {N}onlinear constraints}, Journal of the Society
  for Industrial and Applied Mathematics, 9 (1961), pp.~514--532,
  \url{http://www.jstor.org/stable/2098878}.

\bibitem{Shalit:2012}
{\sc U.~Shalit, D.~Weinshall, and G.~Chechik}, {\em Online learning in the
  embedded manifold of low-rank matrices}, Journal of Machine Learning
  Research, 13 (2012), pp.~429--458.

\bibitem{Simoncini:2007}
{\sc V.~Simoncini}, {\em {A New Iterative Method for Solving Large-Scale
  Lyapunov Matrix Equations}}, SIAM Journal on Scientific Computing, 29 (2007),
  pp.~1268--1288, \url{https://doi.org/10.1137/06066120X}.

\bibitem{Simoncini:2016}
{\sc V.~Simoncini}, {\em Computational methods for linear matrix equations},
  SIAM Review, 58 (2016), pp.~377--441,
  \url{https://doi.org/10.1137/130912839}.

\bibitem{Steinlechner:2016aa}
{\sc M.~Steinlechner}, {\em Riemannian {Optimization} for {High}-{Dimensional}
  {Tensor} {Completion}}, SIAM Journal on Scientific Computing, 38 (2016),
  pp.~S461--S484.

\bibitem{Sutti:2020}
{\sc M.~Sutti and B.~Vandereycken}, {\em {RMGLS: A MATLAB algorithm for
  Riemannian multilevel optimization}}.
\newblock Available online, June 2020,
  \url{https://doi.org/10.26037/yareta:zara3a5aivcsfk6uhq4oovjxhe}.

\bibitem{Toint:2009}
{\sc P.~L. Toint, D.~Tomanos, and M.~Weber-Mendonça}, {\em A multilevel
  algorithm for solving the trust-region subproblem}, Optimization Methods and
  Software, 24 (2009), pp.~299--311,
  \url{https://doi.org/10.1080/10556780802571467}.

\bibitem{Trottenberg:2000}
{\sc U.~Trottenberg, C.~Oosterlee, and A.~Schuller}, {\em Multigrid}, Elsevier
  Science, 2000, \url{https://books.google.ch/books?id=9ysyNPZoR24C}.

\bibitem{Uschmajew:2015}
{\sc A.~{Uschmajew} and B.~{Vandereycken}}, {\em Greedy rank updates combined
  with {R}iemannian descent methods for low-rank optimization}, in 2015
  International Conference on Sampling Theory and Applications (SampTA), May
  2015, pp.~420--424, \url{https://doi.org/10.1109/SAMPTA.2015.7148925}.

\bibitem{UschmajewV:2019}
{\sc A.~Uschmajew and B.~Vandereycken}, {\em Geometric Methods on Low-Rank
  Matrix and Tensor Manifolds}, Springer International Publishing, Cham, 2020,
  ch.~9, pp.~261--313, \url{https://doi.org/10.1007/978-3-030-31351-7_9}.

\bibitem{Vandereycken:2013}
{\sc B.~Vandereycken}, {\em {Low-Rank Matrix Completion by Riemannian
  Optimization}}, SIAM Journal on Optimization, 23 (2013), pp.~1214--1236,
  \url{https://doi.org/10.1137/110845768}.

\bibitem{VandereyckenV_2010}
{\sc B.~Vandereycken and S.~Vandewalle}, {\em A {R}iemannian optimization
  approach for computing low-rank solutions of {L}yapunov equations}, SIAM
  Journal on Matrix Analysis and Applications, 31 (2010), pp.~2553--2579,
  \url{https://doi.org/10.1137/090764566}.

\bibitem{Wen:2009}
{\sc Z.~Wen and D.~Goldfarb}, {\em A line search multigrid method for
  large-scale nonlinear optimization}, SIAM Journal on Optimization, 20 (2009),
  pp.~1478--1503, \url{https://doi.org/10.1137/08071524X}.

\end{thebibliography}
\end{document}